\documentclass[12pt]{article}
\usepackage{amssymb}
\usepackage{amsfonts,amsmath,eucal,amssymb}

\def \beq {\begin{eqnarray}}
\def \eeq {\end{eqnarray}}
\def \beqn {\begin{eqnarray*}}
\def \eeqn {\end{eqnarray*}}

\newcommand{\halmos}{\rule{1ex}{1.4ex}}

\newcounter{for}[section]

\newtheorem{itlemma}{Lemma}[section]
\newtheorem{itproposition}[itlemma]{Proposition}

\newtheorem{theorem}[itlemma]{Theorem}
\newtheorem{itcorollary}[itlemma]{Corollary}
\newtheorem{itremark}[itlemma]{Remark}
\newtheorem{itremarks}[itlemma]{Remarks}
\newtheorem{itdefinition}[itlemma]{Definition}
\newtheorem{itexample}[itlemma]{Example}

\newenvironment{lemma}{\begin{itlemma}}{\end{itlemma}}
\newenvironment{remark}{\begin{itremark}\rm}{\end{itremark}}
\newenvironment{remarks}{\begin{itremarks} \rm}{\end{itremarks}}
\newenvironment{corollary}{\begin{itcorollary}}{\end{itcorollary}}
\newenvironment{proposition}{\begin{itproposition}}{\end{itproposition}}
\newenvironment{definition}{\begin{itdefinition}\rm}{\end{itdefinition}}
\newenvironment{example}{\begin{itexample}\rm}{\end{itexample}}
\newenvironment{proof}{\noindent {\em Proof}.\ \
}{\hspace*{\fill}$\halmos$\medskip}
\newcommand{\be}[1]{\addtocounter{for}{1} \begin{equation}\label{#1}}
\newcommand{\ee}{\end{equation}}
\newcommand{\bl}[1]{\begin{lemma}\label{#1}}
\newcommand{\br}[1]{\begin{remark}\label{#1}}
\newcommand{\brs}[1]{\begin{remarks}\label{#1}}
\newcommand{\bt}[1]{\begin{theorem}\label{#1}}
\newcommand{\bd}[1]{\begin{definition}\label{#1}}
\newcommand{\bp}[1]{\begin{proposition}\label{#1}}
\newcommand{\bfact}[1]{\begin{fact}\label{#1}}
\newcommand{\bc}[1]{\begin{corollary}\label{#1}}
\newcommand{\bex}[1]{\begin{example}\label{#1}}
\newcommand{\ec}{\end{corollary}}
\newcommand{\efact}{\end{fact}}
\newcommand{\eex}{\end{example}}
\newcommand{\el}{\end{lemma}}
\newcommand{\er}{\end{remark}}
\newcommand{\ers}{\end{remarks}}
\newcommand{\et}{\end{theorem}}
\newcommand{\ed}{\end{definition}}
\newcommand{\ep}{\end{proposition}}
\newcommand{\epr}{\end{proof}}
\newcommand{\bpr}{\begin{proof}}
\newcommand{\bcl}[1]{\begin{claim}\label{#1}}
\newcommand{\ecl}{\end{claim}}

\newcommand{\ecs}{\end{corollary}}
\newcommand{\eers}{\end{exercise}}
\newcommand{\eexs}{\end{example}}
\newcommand{\eems}{\end{example}}
\newcommand{\els}{\end{lemma}}
\newcommand{\eles}{\end{lemmaex}}
\newcommand{\ets}{\end{theorem}}
\newcommand{\eds}{\end{definition}}
\newcommand{\eps}{\end{proposition}}

\newcommand{\bi}{\begin{itemize}}
\newcommand{\ei}{\end{itemize}}
\newcommand{\ben}{\begin{enumerate}}
\newcommand{\een}{\end{enumerate}}

\def\vbar{\mathchoice{\vrule height6.3ptdepth-.5ptwidth.8pt\kern-.8pt}
   {\vrule height6.3ptdepth-.5ptwidth.8pt\kern-.8pt}
   {\vrule height4.1ptdepth-.35ptwidth.6pt\kern-.6pt}
   {\vrule height3.1ptdepth-.25ptwidth.5pt\kern-.5pt}}
\def\fudge{\mathchoice{}{}{\mkern.5mu}{\mkern.8mu}}
\def\bbc#1#2{{\rm \mkern#2mu\vbar\mkern-#2mu#1}}
\def\bbb#1{{\rm I\mkern-3.5mu #1}}
\def\bba#1#2{{\rm #1\mkern-#2mu\fudge #1}}
\def\bb#1{{\count4=`#1 \advance\count4by-64 \ifcase\count4\or\bba A{11.5}\or
   \bbb B\or\bbc C{5}\or\bbb D\or\bbb E\or\bbb F \or\bbc G{5}\or\bbb H\or
   \bbb I\or\bbc J{3}\or\bbb K\or\bbb L \or\bbb M\or\bbb N\or\bbc O{5} \or
   \bbb P\or\bbc Q{5}\or\bbb R\or\bbc S{4.2}\or\bba T{10.5}\or\bbc U{5}\or
   \bba V{12}\or\bba W{16.5}\or\bba X{11}\or\bba Y{11.7}\or\bba Z{7.5}\fi}}

\def \qed {{\hspace*{\fill}$\halmos$\medskip}}

\def \Z {{\mathbb Z}}
\def \R {{\mathbb R}}

\def \N {{\mathbb N}}

\def \ra {\rightarrow }

\def \a {\alpha}
\def \O{\Omega}
\def \s {\sigma}
\def \T {{\cal{T}}}

\def\g{\gamma}
\def\d{\delta}
\def\e{\varepsilon}

\def\b{\beta}
\def\L{\Lambda}
\def\l{\lambda}
\def\E{{\cal{E}}}

\numberwithin{equation}{section}

\newtheorem{teo}{Theorem}[section]
\newtheorem{de}[teo]{Definition}

\DeclareMathOperator{\diam}{diam}

\DeclareMathOperator{\gap}{Gap}

\DeclareMathOperator{\ent}{Ent}
\DeclareMathOperator{\Range}{Range}

\newcommand{\ind}{\mathbf{1}}

\renewcommand{\natural}{\mathbb{N}}

\newcommand{\ie}{\emph{i.e.\ }}

\newcommand{\eg}{\emph{e.g.\ }}

\DeclareMathSymbol{\leqslant}{\mathalpha}{AMSa}{"36} 
\DeclareMathSymbol{\geqslant}{\mathalpha}{AMSa}{"3E} 
\DeclareMathSymbol{\eset}{\mathalpha}{AMSb}{"3F}     
\renewcommand{\leq}{\;\leqslant\;}                   
\renewcommand{\geq}{\;\geqslant\;}                   

\newcommand{\grad}{\nabla} 

\newcommand{\la}{\label}

\def\1{\ifmmode {1\hskip -3pt \rm{I}} \else {\hbox {$1\hskip -3pt \rm{I}$}}\fi}

\newcommand{\cE}{\ensuremath{\mathcal E}}

\newcommand{\cL}{\ensuremath{\mathcal L}} 
\newcommand{\cM}{\ensuremath{\mathcal M}}

\newcommand{\bbN}{{\ensuremath{\mathbb N}} }

\newcommand{\bbR}{{\ensuremath{\mathbb R}} }

\newcommand{\bbZ}{{\ensuremath{\mathbb Z}} }

\let\a=\alpha \let\b=\beta   \let\d=\delta  \let\e=\varepsilon
 \let\g=\gamma       \let\l=\lambda
            
\let\r=\rho  \let\s=\sigma    
  
     \let\L=\Lambda 
\let\O=\Omega      

\def\\{\hfill\break}

\def\tthsp{\kern .083333 em}

\def\?{\mskip -10mu}
\def\indbox#1{\hbox to \parindent{\hfil\ #1\hfil} }

\newfam\msafam
\newfam\msbfam
\newfam\eufmfam
\def\hexnumber#1{%
  \ifcase#1 0\or 1\or 2\or 3\or 4\or 5\or 6\or 7\or 8\or
  9\or A\or B\or C\or D\or E\or F\fi}
\font\tenmsa=msam10
\font\sevenmsa=msam7
\font\fivemsa=msam5
\textfont\msafam=\tenmsa
\scriptfont\msafam=\sevenmsa
\scriptscriptfont\msafam=\fivemsa        
\edef\msafamhexnumber{\hexnumber\msafam}%
\mathchardef\restriction"1\msafamhexnumber16
\mathchardef\ssim"0218
\mathchardef\square"0\msafamhexnumber03
\mathchardef\eqd"3\msafamhexnumber2C
\def\QED{\ifhmode\unskip\nobreak\fi\quad
  \ifmmode\square\else$\square$\fi}            
\font\tenmsb=msbm10
\font\sevenmsb=msbm7
\font\fivemsb=msbm5
\textfont\msbfam=\tenmsb
\scriptfont\msbfam=\sevenmsb
\scriptscriptfont\msbfam=\fivemsb

\font\teneufm=eufm10
\font\seveneufm=eufm7
\font\fiveeufm=eufm5
\textfont\eufmfam=\teneufm
\scriptfont\eufmfam=\seveneufm
\scriptscriptfont\eufmfam=\fiveeufm

\def\({\left(}
\def\){\right)}
\let\neper=e
\let\ii=i

\def\ie{\hbox{\it i.e.\ }}

\def\nep#1{ \neper^{#1}}

\let\<=\langle
\let\>=\rangle

\def\diam{\mathop{\rm diam}\nolimits}

\def\gap{\mathop{\rm gap}\nolimits}

\outer\def\nproclaim#1 [#2]#3. #4\par{\medbreak \noindent
   \talato(#2){\bf #1 \Thm[#2]#3.\enspace }%
   {\sl #4\par }\ifdim \lastskip <\medskipamount 
   \removelastskip \penalty 55\medskip \fi}
\def\thmm[#1]{#1}
\def\teo[#1]{#1}
\def\sttilde#1{%
\dimen2=\fontdimen5\textfont0
\setbox0=\hbox{$\mathchar"7E$}
\setbox1=\hbox{$\scriptstyle #1$}
\dimen0=\wd0
\dimen1=\wd1
\advance\dimen1 by -\dimen0
\divide\dimen1 by 2
\vbox{\offinterlineskip%
   \moveright\dimen1 \box0 \kern - \dimen2\box1}
}
\begin{document}

\title{Spectral gap estimates for interacting particle systems via a
  Bochner--type identity}

\author{
Anne-Severine Boudou \\ Dipartimento di Matematica Pura e Applicata \\
Universit\`{a} di Padova \\ Via Belzoni 7, 35131 Padova, Italy \\ 
e-mail: \texttt{boudou@math.unipd.it}\\~
\\ Pietro Caputo \\ Dipartimento di Matematica \\ Universit\`a di Roma 3, L.go S. Murialdo 1, 00146 Roma, Italy \\ e-mail: \texttt{caputo@mat.uniroma3.it}
\\~
\\
Paolo Dai Pra \\
Dipartimento di Matematica Pura e Applicata \\
Universit\`{a} di Padova \\ Via Belzoni 7, 35131 Padova, Italy \\ 
e-mail: \texttt{daipra@math.unipd.it}\\~
\\Gustavo Posta \\
Dipartimento di Matematica \\ Politecnico di Milano \\ Piazza L. Da Vinci 32, 20133 Milano, Italy\\
e-mail: \texttt{gustavo.posta@polimi.it}
}

\date{}

\maketitle

\begin{abstract}
We develop a general technique, based on a Bochner-type identity, to estimate spectral gaps of a class of Markov operator. We apply this technique to various interacting particle systems. In particular, we give a simple and short proof of the diffusive scaling of the spectral gap of the Kawasaki model at high temperature. Similar results are derived for Kawasaki-type dynamics in the lattice without exclusion, and in the continuum. New estimates for Glauber-type dynamics are also obtained.
\end{abstract}

\section{Introduction}

Consider a Markov process $(X_t)_{t \geq 0}$, with values on a measurable space $(S, {\cal S})$, having an invariant measure $\nu$ and whose equilibrium dynamics are time-homogeneous and time-reversible in law. The family of operators 
\[
T_t f(x) := E[f(X_t)|X_0 =x]
\]
form a semigroup of self-adjoint, positivity preserving, contractions on $L^2(\nu)$. Note that $(T_t)_{t \geq 0}$ is well defined, and contractive, in $L^{\infty}(\nu)$ as well, and therefore, by interpolation, on all $L^p(\nu)$ with $2 \leq p \leq +\infty$. One  aim of ergodic theory for Markov process is to understand whether $T_t f$ converges, and in which sense, to the equilibrium average $\nu[f]:= \int f d\nu$ and, if this is the case, to give quantitative estimates on the rate of convergence. One of the main tools in this context is provided by {\em functional inequalities}, in particular {\em Poincar\'e inequality}, \emph{logarithmic-Sobolev inequality} and {\em modified logarithmic-Sobolev inequality}. In order to illustrate the use of these inequalities, assume the semigroup $(T_t)_{t \geq 0}$ has a selfadjoint generator ${\cal L}$ with domain ${{\cal D}}({\cal L})$, as would follow from assuming strong right-continuity. The associated {\em Dirichlet form} is defined on ${{\cal D}}({\cal L}) \times {{\cal D}}({\cal L})$ and is given by 
\[
\E(f,g) := - \nu[f{\cal L}g].
\]
For $f,g \in L^2(\nu)$, let
\[
\nu[f;g] := \nu[fg] - \nu[f]\nu[g]
\]
be the covariance of $f$ and $g$. The inequality
\be{P}
k \, \nu[f;f] \leq \E(f,f) \ \ \mbox{for every } f \in {{\cal D}}({\cal L})
\ee
is called {\em Poincar\'e inequality}. The largest  $k \geq 0$ for which (\ref{P}) holds is the {\em spectral gap} of ${\cal L}$ in $L^2(\nu)$, and we denote it by $\gap({\cal L})$. Indeed, if $\gap({\cal L})>0$, then  (\ref{P}) is equivalent to the fact that $0$ is a simple eigenvalue for ${\cal L}$ (with the constants as eigenvectors), while the remaining part of the spectrum is contained in $(-\infty, -k]$. A straightforward consequence of (\ref{P}) is, therefore,
\[
\|T_t f - \nu[f] \|^2_2 \leq e^{-2kt} \nu[f;f],
\]
for all $f \in L^2(\nu)$, \ie $T_t f$ converges to $\nu[f]$ in $L^2(\nu)$ with exponential rate $\gap({\cal L})$.

Now, let $f \in L^1(\nu)$, $f \geq 0$, and define the {\em entropy}
\[
\ent_{\nu}(f) := \nu[f \log f] - \nu[f]\log \nu[f],
\]
with the conventions $0 \log 0 =0$ and $\ent_{\nu}(f) = +\infty$ if $f \log f \not\in L^1(\nu)$. By Jensen's inequality, it is easily checked that $\ent_{\nu}(f) \geq 0$, and $\ent_{\nu}(f) = 0$ if and only if $f = \mbox{const.}$ $\nu$-a.s.. The inequality
\be{LS}
s \ent(f) \leq \E(\sqrt{f},\sqrt{f})  \ \ \mbox{for every } f \mbox{ such that } \sqrt{f} \in {{\cal D}}({\cal L})
\ee
is called {\em logarithmic-Sobolev inequality}, while
\be{MLS}
\a  \ent(f) \leq \E(f,\log f)  \ \ \mbox{for every } f \mbox{ such that } f, \log f \in {{\cal D}}({\cal L})
\ee
is called {\em modified logarithmic-Sobolev inequality}. These three inequalities are hierarchically ordered in the following sense: if (\ref{LS}) holds with $s>0$, then (\ref{MLS}) holds with $\a \geq s/4$; if (\ref{MLS}) holds with $\a>0$, then (\ref{P}) holds with $k \geq \a/2$. Various consequences of (\ref{LS}) and (\ref{MLS}) in terms of  ergodicity of the semigroup $T_t$ can be obtained (see \eg \cite{Da:Pa:Po}). For instance, under some additional conditions on the domain ${{\cal D}}({\cal L})$, the modified logarithmic-Sobolev inequality is equivalent to the statement
\[
\ent_{\nu}(T_t f) \leq e^{-\a t} \ent_{\nu}(f)
\]
for each $f$ with finite entropy. For diffusion processes, the logarithmic-Sobolev inequality and its modified version coincide, while for Markov processes with discontinuous trajectories the two inequalities are, in general, not equivalent.

The study of functional inequalities for interacting particle systems (\cite{Li}) has been motivated by both theoretical and computational purposes, and has led to the development of a rather sophisticated mathematical technology (\cite{sz1,sz2,Lu:Ya,ma:ol,Ca:Ma:Ro}). The main aim of this paper is to adapt to a class of Markov processes with discontinuous trajectories, including many interesting interacting particle systems, an approach to functional inequalities that goes back to Bochner (\cite{Bo}) and Lichn\'erowicz (\cite{Lic}). This approach was originally developed in the context of Riemannian geometry and allows to obtain lower bounds for  the spectral gap of the Laplacian in Riemannian manifolds. Later
Bakry \& Emery (\cite{Ba:Em}) have used similar ideas in a more general context,  obtaining, in addition to spectral gap estimates,  lower bounds for the best constant in the Logarithmic-Sobolev inequality for diffusion operators. Bakry \& Emery's work has inspired several further developments (e.g. \cite{Ca:St,De:St,He,Le}), in particular concerning diffusion models motivated by statistical mechanics. 

The following proposition is the starting point of the approach we just mentioned.
\bp{bakem}
The spectral gap $\gap({\cal{L}})$ of a Markov generator ${\cal L}$, self-adjoint in $L^2(\nu)$, is equal to the largest constant $k$ such that the inequality
\be{BE}
k\, \E(f,f) \leq \nu\left[ ({\cal L}f)^2\right] 
\ee
holds true for every $ f \in {{\cal D}}({\cal L})$.
\ep
The proof of Proposition \ref{bakem} is  a simple consequence of the spectral Theorem. Indeed, let $(E_{\l})_{\l \geq 0}$ be the spectral projections of the nonnegative, self-adjoint operator $-{\cal{L}}$ in $L^2(\nu)$, and let $k := \gap({\cal{L}})$.
The spectral Theorem yields
\be{sth}
\E(f,f) = \int_{k}^{+\infty} \l \, d(E_{\l}f,f) \ \ \ \  \nu\left[ ({\cal L}f)^2\right] = \int_{k}^{+\infty} \l^2 d(E_{\l}f,f),
\ee
where $(\, \cdot \, , \, \cdot \, )$ denotes here the scalar product in $L^2(\nu)$.
Thus, the inequality 
\eqref{BE} follows from (\ref{sth}) and the obvious fact that $\l^2\geq k \, \l$ on $[k, +\infty)$. \\ In order  to see that $k=\gap({\cal L})$ is the largest constant for which (\ref{BE}) holds for every $ f \in {{\cal D}}({\cal L})$, for a given $\e>0$ we can choose $0 \neq f \in \Range(E_{k+\e}-E_{k^-})$, where $E_{k^-}$ denotes left limit (actually, we choose $0 \neq f \in \Range(E_{\e}-E_{0})$ in the case $k=0$). We have that $f \in  {{\cal D}}({\cal L})$ and, by (\ref{sth}), 
\[
0 < \nu\left[ ({\cal L}f)^2\right] \leq (k+\e)\E(f,f) <  (k+2\e)\E(f,f).
\]
Thus (\ref{BE}) does not hold for $k+2\e$, and the proof of Proposition \ref{bakem} is complete.

In order to obtain explicit estimates for the spectral gap we rewrite the term $\nu\left[ ({\cal L}f)^2\right] $ in a form which can be conveniently compared to the Dirichlet form; in the case of diffusion operators this is realized by the so-called Bochner identity (see \cite{St}, Chapter 6 for a general treatment).
In Section 2 of this paper we prove a version of this identity (Corollary \ref{corollary1}) and
we develop, partly by collecting existing ideas, a general approach to inequality (\ref{BE}) for a very wide class of Markov processes with discontinuous trajectories, including interacting particle systems with a reversible probability measure. We then apply these tools to several models. In Section 3 we prove the diffusive scaling of the spectral gap of the Kawasaki model at sufficiently high temperature. This result goes back to Lu and Yau in \cite{Lu:Ya}; their extremely difficult proof has been made more accessible in \cite{Ca:Ma}, even though it still required a long and technical inductive argument. The statement proved in \cite{Lu:Ya} and \cite{Ca:Ma} is that diffusive scaling of the spectral gap follows from the so-called {\em strong mixing} condition on the associated Gibbs measure, which in turn holds true at sufficiently high temperature (but at any temperature in dimension $d=1$). In this paper we prove the weaker result that diffusive scaling holds at sufficiently high temperature, with no direct connection with mixing properties of the Gibbs measure; although the result is weaker, the proof is quite short and simple. Our approach proves to be very flexible, and has allowed us to give estimates on the spectral gap of other models with conservation of particle number, in particular lattice models with unbounded number of particles (Section 4) and Kawasaki-type dynamics in the continuum (Section 5). For these models spectral gap estimates are not available in the literature. The remaining sections are dedicated to non-conservative models, in particular Glauber dynamics in the lattice with unbounded spin (Section 6), and Glauber dynamics  in the continuum (Section 7). For the models in Section 6, estimates on both spectral gap and the constant in the modified logarithmic-Sobolev inequality were obtained in \cite{Da:Pa:Po}, in the case of uniformly bounded interaction. The method in this paper allows unbounded interaction too. For the models in Section 7, estimates on the spectral gap were obtained first in \cite{Be:Ca:Ce}, via an inductive argument, and then in \cite{Ko:Ly} via the same sort of arguments we use here; our point here is to show that this argument is a special case of a general, and rather powerful, method.

We finally remark that this approach, unlike for diffusion operators, has not yet allowed estimates for the best constant in the logarithmic-Sobolev inequality or its modified version, except for special models (see \cite{Ca:Po}).

\section{General scheme}

In this section we give the formal basis of our method for estimating spectral gaps of a class of Markov dynamics. Suppose $(S, {\cal{S}}, \nu)$ is a probability space. Here $S$ will be interpreted as the state space for the dynamics, and $\nu$ a corresponding invariant probability measure. Let $G$ be a set of measurable transformations from $S$ to $S$, and $ {\cal{G}}$ be a $\s$-field of subsets of $G$. To each $\eta \in S$ we associate a positive $\s$-finite measure $c(\eta,d\g)$ on $(G,{\cal{G}})$ in such a way that for every $
\varphi:G \ra [0,+\infty]$ measurable, the map $\eta \mapsto \int \varphi(\g)c(\eta,d\g)$ is measurable. 
\noindent
In this paper we deal with Markovian dynamics on $S$ whose infinitesimal generator ${\cal L}$ is a well defined linear (possibly unbounded, with dense domain ${\cal D}({\cal L})$) operator on $L^2(\nu)$, given by, for $f \in {\cal D}({\cal L})$
\be{generator}
{\cal L}f(\eta) = \int_G \nabla_{\g}f(\eta) c(\eta,d\g),
\ee
where $\nabla_{\g}f = f \circ \g - f$. This class of operator includes generators of Markov chains with finite or countable state space, as well as interacting particle systems, as defined in Chapter 1 of \cite{Li}. \\ In what follows, $\nu_c$ denotes the positive measure on $S \times G$ given by $\nu_c(d\eta,d\g) := \nu(d\eta)c(\eta,d\g)$. We make the following additional assumption on the generator ${\cal{L}}$.

\vspace{0.5cm}
\noindent
{\bf \em (Rev)} For every $\g \in G$ there is a unique $\g^{-1} \in G$ such that the equality $\g^{-1}(\g(\eta)) = \eta$ holds $\nu_c$-a.s.. Moreover, for every $\Psi \in L^1(\nu_c)$, 
\be{db}
\int \Psi(\eta,\g) c(\eta,d\g) \nu(d\eta) = \int \Psi(\g(\eta),\g^{-1}) c(\eta,d\g) \nu(d\eta).
\ee

\vspace{0.5cm}
Note that assumption {\em (Rev)} implies that ${\cal L}$ is symmetric in $L^2(\nu)$, \ie
\be{sym}
\int f(\eta)g(\g(\eta))c(\eta,d\g)\nu(d\eta) = \int f(\g(\eta))g(\eta)c(\eta,d\g)\nu(d\eta) .
\ee
Thus, {\em (Rev)} is a reversibility condition, and (\ref{db}) is the usual {\em detailed balance} condition written in this general context. 
\\ Note that, under  {\em (Rev)}, for $f,g \in {\cal D}({\cal L})$,
\be{dir}
\E(f,g) := - \nu\left(f{\cal L} g\right) = \frac{1}{2}\nu\left[ \int c(\eta,d\g) \nabla_{\g}f(\eta) \nabla_{\g}g(\eta)\right].
\ee

The method we present in this section is based on the possibility of constructing a positive measure $R$ on $S \times G \times G$ having the following properties.
\bi
\item[{\bf (A1)}]
There is a core ${\cal{C}}$ of ${\cal D}({\cal L})$ such that for each $f \in {\cal{C}}$, the function $(\eta,\g,\d) \mapsto \nabla_{\g}f(\eta) \nabla_{\d}f(\eta)$ belongs to $L^1(R)$.

\item[{\bf (A2)}]
The equality
\[
\g(\d(\eta) ) = \d(\g(\eta))
\]
holds $R$-almost everywhere.

\item[{\bf (A3)}]
Define $\Theta F(\eta,\g,\d) := F(\eta,\d,\g)$, then for any $F\in L^1(R)$.
\begin{displaymath}
  \int \Theta F dR=
  \int  F\;dR .
\end{displaymath}

\item[{\bf (A4)}]

Define $\T F(\eta,\g,\d) := F(\g(\eta),\g^{-1},\d)$, then for any $F\in L^1(R)$.
\begin{displaymath}
  \int \T F\;dR=
  \int  F dR .
\end{displaymath}
\ei

The basic computation is given in the following Lemma.
\bl{lemma1}
For all $f \in {\cal{C}}$
\[
\int \left[ \nabla_{\g}\nabla_{\d}f(\eta) \right]^2 dR = 4\int \nabla_{\g}f(\eta) \nabla_{\d}f(\eta) dR.
\]
\el
\bpr
First, by (A2),
\[
\int \left[ \nabla_{\g}\nabla_{\d}f(\eta) \right]^2 dR = \int  \nabla_{\g}\nabla_{\d}f(\eta)
 \nabla_{\d}\nabla_{\g}f(\eta)dR.
\]
We now write
\begin{multline}
 \nabla_{\g}\nabla_{\d}f(\eta)
 \nabla_{\d}\nabla_{\g}f(\eta) = \\
 \nabla_{\d}f(\g(\eta))\nabla_{\g}f(\d(\eta)) - \nabla_{\d}f(\g(\eta))\nabla_{\g}f(\eta)
 - \nabla_{\d}f(\eta)\nabla_{\g}f(\d(\eta)) + \nabla_{\d}f(\eta)\nabla_{\g}f(\eta). \label{e1}
\end{multline}
We show that each one of the four summands in the r.h.s. of (\ref{e1}) is in $L^1(R)$, and its integral with respect to $R$ equals
\[
\int \nabla_{\d}f(\eta)\nabla_{\g}f(\eta) dR.
\]
From this fact the conclusion follows.
By assumption (A1), for the fourth summand there is nothing to prove. Moreover, using assumption (A4) in the first equality,
\begin{displaymath}
\int \nabla_{\d}f(\eta) \nabla_{\g} f(\eta)dR  = \int \nabla_{\d}  f(\g(\eta)) \nabla_{\g^{-1}}f(\g(\eta))dR  = 
- \int \nabla_{\d} f(\g(\eta)) \nabla_{\g}f(\eta)dR ,
\end{displaymath}
that takes care of the second summand in (\ref{e1}). The integral of the third summand equals the one of the second by assumption (A3). Finally, using first (A4), then (A3), (A4) again and (A2),
\begin{multline*}
\int \nabla_{\d}f(\eta) \nabla_{\g} f(\eta)dR = \int  \nabla_{\d}f(\g(\eta))  \nabla_{\g^{-1}}f(\g(\eta)) dR 
= -\int  \nabla_{\g}f(\d(\eta))  \nabla_{\d}f(\eta) dR = \\
- \int \nabla_{\g^{-1}}f(\d(\g(\eta))) \nabla_{\d}f(\g(\eta)) dR  = \int \nabla_{\g}f(\d(\eta)) \nabla_{\d}f(\g(\eta)) dR .
\end{multline*}
\epr
For a easier reading of the consequences of Lemma \ref{lemma1}, we make the following further assumption.
\bi
\item[{\bf (A5)}]
The measure $R$ is absolutely continuous with respect to the measure $\nu(d\eta)c(\eta,d\g)c(\eta,d\d)$. We denote by $r(\eta,\g,\d)$ the corresponding Radon-Nikodym derivative.
\ei
\bc{corollary1}
For all $f \in {\cal{C}}$
\[
\nu\left[({\cal L}f)^2\right] - \frac{1}{4}\int \left[ \nabla_{\g}\nabla_{\d}f(\eta) \right]^2 dR = 
\int \nu(d\eta)c(\eta,d\g)
c(\eta,d\d)[1-r(\eta,\g,\d)] \nabla_{\g}f(\eta)\nabla_{\d}f(\eta).
\]
\ec
\bpr
It is enough to observe that
\[
\nu\left[({\cal L}f)^2\right] = \int \nu(d\eta)c(\eta,d\g)
c(\eta,d\d) \nabla_{\g}f(\eta)\nabla_{\d}f(\eta),
\]
and apply Lemma \ref{lemma1}.
\epr
Therefore, by Proposition \ref{bakem}, we get the following result.
\bc{corollary2}
If, for all $f$
\[
\int \nu(d\eta)c(\eta,d\g)
c(\eta,d\d)[1-r(\eta,\g,\d)] \nabla_{\g}f(\eta)\nabla_{\d}f(\eta) \geq k \int  \nu(d\eta)c(\eta,d\g)\left[
\nabla_{\g}f(\eta)\right]^2,
\]
then $\gap({\cal L}) \geq 2k$.
\ec
The idea is now to compare, {\em pointwise} in $\eta$, the quadratic forms in $\nabla f$
\be{e2}
\int c(\eta,d\g)
c(\eta,d\d)[1-r(\eta,\g,\d)] \nabla_{\g}f(\eta)\nabla_{\d}f(\eta)
\ee
and 
\be{e3}
\int  c(\eta,d\g)\left[
\nabla_{\g}f(\eta)\right]^2.
\ee
A ``good'' choice for $r(\eta,\g,\d)$ should be when $1-r(\eta,\g,\d)$ is concentrated near the ``diagonal'' $\g = \d$. The following choice works in many examples, including those in Sections 3, 6 and 7 of this paper. In Sections 4 and 5 the $r(\eta,\g,\d)$ given in the Proposition \ref{proposition1} below will need a slight adaptation to the dynamics.
The following additional assumption is needed. 
\bi
\item[{\bf (A6)}]
For $\nu$-almost every $\eta \in S$ and for all $\g \in G$, the measure $c(\g(\eta),d\d)$ is absolutely continuous with respect to the measure $c(\eta,d\d)$.
\ei
\bp{proposition1}
Let us write $G$ in the form $G=J \cup J^{-1}$, where $J \subseteq G$, and $J^{-1} := \{\g : \g^{-1} \in J\}$. $J$ and $J^{-1}$ are not necessarily disjoint. Suppose the reversibility condition (Rev) is satisfied, as well as condition (A6).  Define $r(\eta,\g,\d)$ as follows:
\[
r(\eta,\g,\d) = \left\{
\begin{array}{ll}
\frac{1}{2} \left(1 + \frac{dc(\g(\eta),\cdot)}{dc(\eta,\cdot)}(\d)\right) & \mbox{if } \g \circ \d = \d \circ \g, \ \g,\d \in J \cap J^{-1} \\
 \frac{dc(\g(\eta),\cdot)}{dc(\eta,\cdot)}(\d) & \mbox{if } \g \circ \d = \d \circ \g, \ \g,\d \in J \setminus J^{-1} \mbox{ or } \g,\d \in J^{-1} \setminus J \\
1 & \mbox{if } \g \circ \d = \d \circ \g, \ \left\{ \begin{array}{l} 
        \g \in J \setminus J^{-1}, \  \d \in J^{-1} \setminus J \\ 
         \mbox{or } \g \in J^{-1}\setminus J, \ \d \in J \setminus J^{-1} \end{array} \right. \\
0 & \mbox{otherwise. } 
\end{array} \right.
\]
Then condition (A2) and (A4) are satisfied.
\ep
\bpr
Note that $r(\eta,\g,\d)$ is supported on the set $\{(\eta,\g,\d): \g \circ \d = \d \circ \g\}$, so that (A2) holds easily. To check condition (A4), let $G(\eta,\g,\d)$ be a nonnegative, measurable function. The key fact is given in the following two computations.
\be{e4}
\int \nu(d\eta)c(\eta,d\g)c(\eta,d\d) G(\eta,\g,\d) = \int \nu(d\eta)c(\eta,d\g)c(\g(\eta),d\d) G(\g(\eta),\g^{-1},\d),
\ee
where we have applied (Rev) to the function $\Psi(\eta,\g) := \int c(\eta,d\d) G(\eta,\g,\d)$, and
\be{e5}
\int \nu(d\eta)c(\eta,d\g)c(\g(\eta),d\d) G(\eta,\g,\d) = \int \nu(d\eta)c(\eta,d\g)c(\eta,d\d) G(\g(\eta),\g^{-1},\d),
\ee
where (Rev) has been applied to $\Psi(\eta,\g) := \int c(\g(\eta),d\d) G(\eta,\g,\d)$. \\
Now, let $F(\eta,\g,\d)$ be a nonnegative, measurable function. We have, by (\ref{e4}) and (\ref{e5}),
\begin{multline}
\int_S \int_{\left(J \cap J^{-1}\right)^2} \nu(d\eta)c(\eta,d\g)c(\eta,d\d) r(\eta,\g,\d)F(\eta,\g,\d)\\
 = \frac{1}{2}\int_S\int_{\left(J \cap J^{-1}\right)^2} \nu(d\eta)c(\eta,d\g)c(\eta,d\d) F(\eta,\g,\d)\\
 +\frac{1}{2}\int_S\int_{\left(J \cap J^{-1}\right)^2} \nu(d\eta)c(\eta,d\g)c(\g(\eta),d\d) F(\eta,\g,\d) \\ = 
\frac{1}{2}\int_S\int_{\left(J \cap J^{-1}\right)^2}\nu(d\eta)c(\eta,d\g)c(\g(\eta),d\d) F(\g(\eta),\g^{-1},\d) + \\ \frac{1}{2}\int_S\int_{\left(J \cap J^{-1}\right)^2}\nu(d\eta)c(\eta,d\g)c(\eta,d\d) F(\g(\eta),\g^{-1},\d)\\ = \int_S\int_{\left(J \cap J^{-1}\right)^2}\nu(d\eta)c(\eta,d\g)c(\eta,d\d) r(\eta,\g,\d)F(\g(\eta),\g^{-1},\d) \label{e6}.
\end{multline}
Similarly:
\begin{multline}
\int_S \int_{\left(J \setminus J^{-1}\right)^2} \nu(d\eta)c(\eta,d\g)c(\eta,d\d) r(\eta,\g,\d)F(\eta,\g,\d) \\ =  \int_S \int_{\left(J \setminus J^{-1}\right)^2} \nu(d\eta)c(\eta,d\g)c(\g(\eta),d\d) F(\eta,\g,\d) \\ = \int_S \int_{\left(J^{-1} \setminus J \right) \times \left(J \setminus J^{-1}\right)} \nu(d\eta)c(\eta,d\g)c(\eta,d\d)F(\g(\eta),\g^{-1},\d) \\ = \int_S \int_{\left(J^{-1} \setminus J \right) \times \left(J \setminus J^{-1}\right)}  \nu(d\eta)c(\eta,d\g)c(\eta,d\d) r(\eta,\g,\d)F(\g(\eta),\g^{-1},\d) \label{e7}.
\end{multline}
All other cases are obvious modifications of (\ref{e6}) and (\ref{e7}).
\epr
The integrability assumption (A1) is usually not harmful, the symmetry condition (A3) with the $r(\eta,\g,\d)$ above, depends on the actual choice of the rates $c(\eta,d\g)$.

\br{diagonal}
In some cases a  modification of the $r(\eta,\g,\d)$ given in Proposition \ref{proposition1} is convenient. Consider the set
\[
D := \{(\eta,\g,\d) \in S \times G \times G: \g = \g^{-1} = \d\}.
\]
Note that both $D$ and $D^c$ are stable for the maps $\Theta$ and $\T$. Therefore we can force $r(\eta,\g,\d) \equiv 0$ for $(\eta,\g,\d) \in D$ without modifying the validity of properties (A1)-(A4).
\er

\br{dirichlet}
In the setting above, if $J \cap J^{-1} = \emptyset$,  useful expressions for the Dirichlet form (\ref{dir}) are
\be{dirj}
\E(f,g) = \nu\left[ \int_J c(\eta,d\g) \nabla_{\g}f(\eta) \nabla_{\g}g(\eta)\right] = \nu\left[ \int_{J^{-1}}c(\eta,d\g) \nabla_{\g}f(\eta) \nabla_{\g}g(\eta)\right],
\ee
as is easily checked using (\ref{db}).
\er

\section{The Kawasaki model}

For a given finite $\L \subset \Z^d$, we consider a model with the finite state space $S:=\{0,1\}^{\L}$. Therefore $\eta \in S$ is of the form $(\eta_x)_{x \in \L}$, where $\eta_x \in \{0,1\}$ is the {\em occupation number} at $x \in \L$. The only allowed transitions are the exchanges of the occupation numbers in two distinct sites $x,z \in \L$. If $\g$ is such exchange map, we write $\g = xz$ and
$\g(\eta) = \eta^{xz}$. So we let 
\[
G = \{ xz: x,z \in \L, \, x \neq z \}.
\]
Let $\Phi = (\Phi_A)_{A \subset \subset \Z^d}$ be a {\em summable potential} in $\Z^d$, \ie for all $A$ finite subset of $\Z^d$, $\Phi_A : \{0,1\}^A \ra \R$, and
\[
\|\Phi\| := \sup_{x \in \Z^d} \sum_{A \ni x} \sup_{\eta} |\Phi_A(\eta)| < +\infty.
\]
In this section we impose the following stronger summability condition
\be{summ}
|||\Phi |||  := \sup_{x \in \Z^d} \sum_{A \ni x} |A|\sup_{\eta} |\Phi_A(\eta)| < +\infty.
\ee
Note that we are {\em not} assuming the potential to be translation invariant or of finite range. \\
Now let $\eta \in S$ and $\tau \in \{0,1\}^{\L^c}$. The element $\eta \tau \in \{0,1\}^{\Z^d}$ is then defined by $(\eta \tau)_x = \eta_x$ for $x \in \L$, and $(\eta \tau)_x = \tau_x$ for $x \in \L^c$. The {\em energy} of $\eta \in S$ is defined by
\[
H^{\tau}(\eta) =  \sum_{A: A \cap \L \neq \emptyset} \Phi_A(\eta \tau),
\]
In the sequel, the {\em boundary condition} $\tau$ will be omitted: indeed, all estimates will be uniform in the boundary conditions. \\
In this section we consider the Kawasaki model in the {\em complete graph}, \ie exchanges in the occupation numbers may occur in {\em any} pair of sites $x,z \in \L$. We study the dynamics determined by the following infinitesimal generator:
\be{e8}
{\cal L}f(\eta) = \sum_{xz} c(\eta,xz)\nabla_{xz}f(\eta)
\ee
where the sum ranges over all pairs $x,z \in \L$, and with
\be{ratesk}
c(\eta,xz) = \frac{1}{|\L|} e^{-\frac{1}{2} \b \nabla_{xz}H(\eta)},
\ee
where $\b>0$ is the inverse temperature. This dynamics conserve the number of occupied sites. For every $0 \leq N \leq |\L|$ we consider the set $S_N$ of configurations with $N$ particles, \ie
\[
S_N := \{ \eta \in S: \sum_{x \in \L} \eta_x = N\},
\]
and the probability measure
\[
\nu_\Lambda^N (\eta) = \frac{1}{Z_\Lambda^N} e^{-\b H(\eta)} \ind(\eta\in S_N),
\]
where $Z_\Lambda^N$ is a normalization factor. All $\nu_\Lambda^N$ are invariant for the dynamics, and the {\em detailed balance condition}
\be{detbal}
c(\eta,xz)\nu_\Lambda^N(\eta) = c(\eta^{xz},xz)\nu_\Lambda^N(\eta^{xz})
\ee
is satisfied. 
We state the main result of this section.
\bt{kawasaki}
Consider the system with the generator ${\cal L}$ in (\ref{e8}), with state space $S_N$ and invariant measure $\nu :=\nu_\Lambda^N$. Assume condition (\ref{summ}) holds. 
For every $0 < \l <1$ there exists $\b_{\l}>0$, independent of $\L,\tau,N$,  such that for $\b \leq \b_{\l}$ we have $\gap({\cal L}) \geq \l$.
\et
\br{note}
For $\b = 0$ the model reduces to simple exclusion in the complete graph, whose gap is known to be equal to $1$. Thus the lower bound in Theorem \ref{kawasaki} becomes optimal in the limit $\b \ra 0$.
\er
Usually, rather than the generator in (\ref{e8}), one considers dynamics where only exchanges between nearest neighbors are allowed:
\[
{\cal L}^{n.n.}f(\eta) = \sum_{x\sim z} e^{-\frac{1}{2} \b \nabla_{xz}H(\eta)} \nabla_{xz}f(\eta),
\]
where the sum $\sum_{x\sim z}$ ranges over pairs $x,z \in \L$ with $|x-z|=1$. In the case the potential $\Phi$ is of finite range, \ie $\Phi_A \equiv 0$ up to a finite number of sets $A$, Lemma 4.3 in \cite{Ya} can be used in a standard way to connect the gap of ${\cal L}^{n.n.}$ with that of ${\cal L}$, getting the following result.
\bc{kawasakinn}
Let $\diam(\L)=\max\{|x-z|: x,z \in \L\}$, and assume $\Phi$ is a finite range potential. There exists $\bar{\b}>0$ and a constant $C>0$, both independent of $\L,\tau,N$, such that for every $\b \leq \bar{\b}$ we have $\gap({\cal L}^{n.n.}) \geq C/\diam(\L)^2$.
\ec
In order to prove Theorem \ref{kawasaki}, we use Corollary \ref{corollary2} with the choice of $R$ as in Proposition \ref{proposition1}, with the modification given in Remark \ref{diagonal}. Note that, in this model, each $\g \in G$ coincides with its inverse. So we are forced to choose $J = J^{-1} = G$. Note that two exchanges $xz$ and $yu$ commute if and only if either $xz = yu$ or $\{x,z\} \cap \{y,u\} = \emptyset$. Thus we get
\be{e9}
R(\eta,xz,yu) = \begin{cases} \nu(\eta) c(\eta,xz) \frac{c(\eta,yu) + c(\eta^{xz},yu)}{2} & \mbox{if } \{x,z\} \cap \{y,u\} = \emptyset \\ 0 & \mbox{otherwise.}
\end{cases}
\ee
\bl{kawl1}
For the measure $R$ given in (\ref{e9}), properties (A1)-(A4) are satisfied.
\el
\bpr
Property (A1) is trivial, since both $S$ and $G$ are finite sets. The reversibility condition (Rev) (see (\ref{db})) is a simple consequence of (\ref{detbal}). Properties (A2) and (A4) are guaranteed by Proposition \ref{proposition1}. The symmetry property (A3) follows from the fact that the quantity $c(\eta,xz)c(\eta^{xz},yu)$ is symmetric in $xz,yu$, as one checks using (\ref{ratesk}).
\epr
{\bf Proof of Theorem \ref{kawasaki}}. By Corollary \ref{corollary2} it is enough to check that
\be{cor}
\int \nu(d\eta)c(\eta,d\g)
c(\eta,d\d)[1-r(\eta,\g,\d)] \nabla_{\g}f(\eta)\nabla_{\d}f(\eta) \geq k(\b) \int  \nu(d\eta)c(\eta,d\g)\left[
\nabla_{\g}f(\eta)\right]^2,
\ee
where $k(\b) \ra \frac{1}{2}$ as $\b \ra 0$. We first note that
\begin{multline}
\int \nu(d\eta)c(\eta,d\g)
c(\eta,d\d)[1-r(\eta,\g,\d)] \nabla_{\g}f(\eta)\nabla_{\d}f(\eta) =\\ 
\sum_{xz, zu} \nu\left[c(\eta,xz)c(\eta,zu) \nabla_{xz}f(\eta)\nabla_{zu}f(\eta)\right] 
 \\
+ 
\frac{1}{2} \sum_{xz,yu:\{x,z\} \cap \{y,u\} = \emptyset} \nu\left[ c(\eta,xz)c(\eta,yu)\left(1-\frac{c(\eta^{xz},yu)}{c(\eta,yu)} \right)\nabla_{xz}f(\eta)\nabla_{yu}f(\eta)\right] . \label{a}
\end{multline}
It is useful to keep in mind that in (\ref{a}) we are summing over elements of $G$, so that, for instance, $xz$ and $zx$ are equal. In particular, the sum $\sum_{xz, zu}$ runs over pairs in $G \times G$ whose corresponding exchanges involve {\em at least} one common point. \\
We estimate the two summands in the r.h.s. of (\ref{a}) separately. We begin by showing the following identity:
\begin{multline}
\sum_{xz,zu} \nu\left[c(\eta,zu) \nabla_{xz}f(\eta)\nabla_{zu}f(\eta)\right] = \frac{|\L|}{2} \sum_{zu} \nu\left[c(\eta,zu) \left(\nabla_{zu}f(\eta)\right)^2\right] \\ = \frac{|\L|}{2} \int  \nu(d\eta)c(\eta,d\g)\left[
\nabla_{\g}f(\eta)\right]^2. \label{b}
\end{multline}
The second equality in (\ref{b}) is obvious. For the first, observe that
\begin{multline}
\sum_{xz,zu} \nu\left[c(\eta,zu) \nabla_{xz}f(\eta)\nabla_{zu}f(\eta)\right] \\ = \sum_{xz \neq zu} \nu\left[c(\eta,zu) \nabla_{xz}f(\eta)\nabla_{zu}f(\eta)\right] + \sum_{zu}\nu\left[c(\eta,zu) \left(\nabla_{zu}f(\eta)\right)^2\right].
\label{b1}
\end{multline}
By (Rev) we have
\begin{multline}
\sum_{xz \neq zu} \nu\left[c(\eta,zu) \nabla_{xz}f(\eta)\nabla_{zu}f(\eta)\right] = \sum_{xz \neq zu}  \nu\left[c(\eta,zu) \nabla_{xz}f(\eta)\nabla_{zu}f(\eta)1_{\{\eta_x \neq \eta_z\} \cap \{\eta_u \neq \eta_z\}}\right]
\\ = - \sum_{xz \neq zu}  \nu\left[c(\eta,zu) \nabla_{xz}f(\eta^{zu})\nabla_{zu}f(\eta)1_{\{\eta_x \neq \eta_u\} \cap \{\eta_u \neq \eta_z\}}\right]  \\ = - \sum_{xz \neq zu}  \nu\left[c(\eta,zu) \nabla_{xu}f(\eta)\nabla_{zu}f(\eta)\right] + \sum_{xz \neq zu}  \nu\left[1_{\{\eta_x \neq \eta_u\} }c(\eta,zu) \left(\nabla_{zu}f(\eta)\right)^2\right], \label{c}
\end{multline}
where we have used the fact that in the set $\{\eta_x \neq \eta_u\} \cap \{\eta_u \neq \eta_z\} $ the identity $ \nabla_{xz}f(\eta^{zu})  = \nabla_{xu}f(\eta) - \nabla_{zu}f(\eta)$ holds. Now note that in the sum $\sum_{xz \neq zu}  \nu\left[c(\eta,zu) \nabla_{xu}f(\eta)\nabla_{zu}f(\eta)\right] $ the condition $xz \neq zu$, i.e. $x \neq u$ does not play any role since, if $x=u$, then $\nabla_{xu}f(\eta) \equiv 0$. However, in the same sum, $xz$ is an element of $G$, which means $x \neq z$ or, equivalently, that $xu$ and $zu$ have {\em exactly} one common point. It follows that
\be{c1}
\sum_{xz \neq zu}  \nu\left[c(\eta,zu) \nabla_{xu}f(\eta)\nabla_{zu}f(\eta)\right] = \sum_{xz \neq zu} \nu\left[c(\eta,zu) \nabla_{xz}f(\eta)\nabla_{zu}f(\eta)\right].
\ee
Moreover, by (Rev), 
\[
\sum_{xz \neq zu}  \nu\left[1_{\{\eta_x \neq \eta_u\} }c(\eta,zu) \left(\nabla_{zu}f(\eta)\right)^2\right] = \sum_{xz \neq zu}  \nu\left[1_{\{\eta_x = \eta_u\} }c(\eta,zu) \left(\nabla_{zu}f(\eta)\right)^2\right],
\]
so that
\begin{multline}\label{c2}
\sum_{xz \neq zu}  \nu\left[1_{\{\eta_x \neq \eta_u\} }c(\eta,zu) \left(\nabla_{zu}f(\eta)\right)^2\right] = \frac{1}{2}\sum_{xz \neq zu}  \nu\left[c(\eta,zu) \left(\nabla_{zu}f(\eta)\right)^2\right] \\
= (|\L| - 2) \sum_{zu} \nu\left[c(\eta,zu) \left(\nabla_{zu}f(\eta)\right)^2\right],
\end{multline}
where we have used the fact that, for a fixed $zu \in G$, the number of elements of $G$ with {\em exactly} one point in common with $zu$ is $2(|\L|-2)$.
Thus, inserting (\ref{c2}) and (\ref{c1}) in (\ref{c}) we get
\be{c3}
\sum_{xz \neq zu} \nu\left[c(\eta,zu) \nabla_{xz}f(\eta)\nabla_{zu}f(\eta)\right] = \left(\frac{|\L|}{2} -1\right)\sum_{zu} \nu\left[c(\eta,zu) \left(\nabla_{zu}f(\eta)\right)^2\right],
\ee
that, inserted in (\ref{b1}) yields (\ref{b}).

Now, let $\e:= \b\|\Phi\|$.  Since
\[
|\b \nabla_{xz}H(\eta)| = \b \left|\sum_{A \cap \{x,z\} \neq \emptyset} \nabla_{xz} \Phi_A(\eta) \right| \leq 2 \e,
\]
we have
\be{d}
\frac{1}{|\L|} e^{-\e} \leq c(\eta,xz) \leq \frac{1}{|\L|} e^{\e}.
\ee
Thus, by (\ref{b}), (\ref{d}) and Schwarz inequality,
\begin{multline}
\sum_{xz,zu} \nu\left[c(\eta,xz)c(\eta,zu) \nabla_{xz}f(\eta)\nabla_{zu}f(\eta)\right] \\ = \frac{1}{|\L|}\sum_{xz,zu} \nu\left[c(\eta,zu) \nabla_{xz}f(\eta)\nabla_{zu}f(\eta)\right] + \sum_{xz,zu} \nu\left[\left(c(\eta,xz)-\frac{1}{|\L|}\right)c(\eta,zu) \nabla_{xz}f(\eta)\nabla_{zu}f(\eta)\right]  \\ = \frac{1}{2} \sum_{zu} \nu\left[c(\eta,zu) \left(\nabla_{zu}f(\eta)\right)^2\right]  + \sum_{xz,zu} \nu\left[\left(c(\eta,xz)-\frac{1}{|\L|}\right)c(\eta,zu) \nabla_{xz}f(\eta)\nabla_{zu}f(\eta)\right]  \\    \geq \frac{1}{2} \sum_{zu} \nu\left[c(\eta,zu) \left(\nabla_{zu}f(\eta)\right)^2\right]  \\ - \sum_{xz,zu} \nu\left[\left|c(\eta,xz)-\frac{1}{|\L|}\right|\sqrt{\frac{c(\eta,zu)}{c(\eta,xz)} }\sqrt{c(\eta,zu)c(\eta,xz)}\left|\nabla_{xz}f(\eta)\nabla_{zu}f(\eta)\right|\right]  \\ \geq \frac{1}{2} \sum_{zu} \nu\left[c(\eta,zu) \left(\nabla_{zu}f(\eta)\right)^2\right] - e^{\e}\left(e^{\e}-1\right) \sum_{zu} \nu\left[c(\eta,zu) \left(\nabla_{zu}f(\eta)\right)^2\right] \\ = \left[\frac{1}{2} - e^{\e}\left(e^{\e}-1\right)\right] \sum_{zu} \nu\left[c(\eta,zu) \left(\nabla_{zu}f(\eta)\right)^2\right].
\label{e}
\end{multline}
This takes care of the first summand in the r.h.s. of (\ref{a}). We now deal with the second summand in the r.h.s. of  (\ref{a}). First we note that
\[
\frac{c(\eta^{xz},yu)}{c(\eta,yu)} = \exp\left[ - \frac{1}{2} \sum_{\stackrel{A \cap \{y,u\} \neq \emptyset}{\scriptscriptstyle A \cap \{x,z\} \neq \emptyset}} \nabla_{xz} \nabla_{yu} \Phi_A(\eta)\right].
\]
Thus, using the inequality $\left|e^x -1\right| \leq |x|e^{|x|}$, we get
\be{f}
\left|1-\frac{c(\eta^{xz},yu)}{c(\eta,yu)} \right| \leq \frac{1}{2} \b \sum_{\stackrel{A \cap \{y,u\} \neq \emptyset}{\scriptscriptstyle A \cap \{x,z\} \neq \emptyset}} \nabla_{xz} \nabla_{yu} \left|\Phi_A(\eta)\right| e^{4\e}.
\ee
On the other hand
\[
\sum_{xz: A \cap \{x,z\} \neq \emptyset} \nabla_{xz} \nabla_{yu} \left|\Phi_A(\eta)\right| \leq 8 |\L||A| \sup_{\eta}\left|\Phi_A(\eta)\right|.
\]
Thus, by (\ref{f})
\be{g}
\sum_{xz: \{x,z\} \cap \{y,u\} = \emptyset} \left|1-\frac{c(\eta^{xz},yu)}{c(\eta,yu)} \right| \leq 4\b |\L| e^{4\e} \sum_{A:A \cap \{y,u\} \neq \emptyset} |A| \sup_{\eta}\left|\Phi_A(\eta)\right| \leq 8 \b |\L|e^{4\e} \||\Phi\||.
\ee
Therefore, by (\ref{g}), and using Schwarz inequality as in (\ref{g}),
\begin{multline}
\left| \sum_{xz,yu:\{x,z\} \cap \{y,u\} = \emptyset} \nu\left[ c(\eta,xz)c(\eta,yu)\left(1-\frac{c(\eta^{xz},yu)}{c(\eta,yu)} \right)\nabla_{xz}f(\eta)\nabla_{yu}f(\eta)\right]  \right| \\ \leq 4 \b e^{5\e} \||\Phi\|| \sum_{u,z} \nu\left[c(\eta,zu) \left(\nabla_{zu}f(\eta)\right)^2\right].
\label{h}
\end{multline}
Finally, by (\ref{a}), (\ref{e}) and (\ref{h}), we get
\begin{multline*}
\int \nu(d\eta)c(\eta,d\g)
c(\eta,d\d)[1-r(\eta,\g,\d)] \nabla_{\g}f(\eta)\nabla_{\d}f(\eta) \\ \geq \left[\frac{1}{2} - e^{\e}\left(e^{\e}-1\right)- 4 \b e^{5\e} \||\Phi\||  \right] \sum_{u,z} \nu\left[c(\eta,zu) \left(\nabla_{zu}f(\eta)\right)^2\right],
\end{multline*}
from which (\ref{cor}) follows.
\qed

\section{Random walks on the complete graph}
Random walks on the complete graph interacting via a zero--range
potential were considered in \cite{Ca:Po}. It was shown that the
spectral gap of the process is positive as soon as a uniform
log--concavity assumption is satisfied. Here we consider
the case where we add a non--zero--range 
interaction to the system. It turns out that the general
method described in the previous sections gives interesting 
conclusions for a wide class of models. 


The reference model is the zero--range process 
obtained as follows. We denote by $V_n$ the
the set of $n$ labeled vertexes and  
consider $N$ random walks on the complete graph
over $V_n$, \ie\ a process of $N$ particles taking 
jumps between any pair of vertexes of $V_n$. 
The state space is  
$$
S_N:=\left\{\eta:V_n\to\N\;\,\text{such that}\,\; \sum_{x\in V_n}\eta_x=N\right\}\,,$$
with $\eta_x$ representing the number of particles at vertex $x$.
At each vertex $x\in V_n$ we associate 
a rate function 
$g_x:\bbN\to \bbR$ such that $g_x(0)=0$, 
and 
\begin{equation}
\min_{x\in V_n} \inf_{k\geq 1} g_x(k)\geq 1\,.
\label{d0}
\end{equation}
The choice in (\ref{d0}) is purely conventional and 
any positive constant instead of $1$ can be accepted 
(this amounts to 
a trivial time rescaling).
A particle is moved from $x$ to a 
uniformly chosen vertex $z\in V_n$ 
with rate $g_x(\eta):= g_x(\eta_x)$ and the 
Markov generator can be written as 
\begin{equation}
\cL f = \frac1n\sum_{x,z} g_x\grad_{xz}f\,,
\label{gen00}
\end{equation}
with the sum extending over all $x,z\in V_n$. Here 
$\grad_{xz}f$ stands for the gradient $f^{xz} -f$, 
with $f^{xz}(\eta) = f(\eta^{xz})$, $\eta^{xz}$ being 
the configuration in which a 
particle has been moved from $x$ to $z$, \ie\ 
\begin{displaymath}
  (\eta^{xz})_y:=
  \begin{cases}
    \eta_x-1 & \text{if $\eta_x>0$ and $y=x$}\\
    \eta_z+1 & \text{if $\eta_x>0$ and $y=z$}\\
    \eta_y   & \text{otherwise.}
  \end{cases}
\end{displaymath}
In this way $\eta^{xy}=\eta$ if $\eta_x=0$. We also agree that
$\eta^{xz}=\eta$
when $x=z$.  When $N=1$ we have a random walk on the 
(weighted) complete graph.
For $N\geq 2$, if the functions $g_x$ were all 
linear, \ie $g_x(n)=\d_x n$ for some constants $\d_x>0$, 
the resulting $N$ random walks would be {\em independent}. Under 
the only assumption (\ref{d0}), however,
in general there is non--trivial interaction.  
The process is reversible w.r.t.\ the probability measure 
$\bar\nu_{V_n}^{N}$ on $S_N$ given by 
\begin{displaymath}
  \bar\nu_{V_n}^N(\eta):=
  \frac{1}{\bar Z_{V_n}^N\prod_x [g_x(\eta_x)!]}\,
\end{displaymath}
where $[g_x(k)!]:=\prod_{\ell=1}^k g_x(\ell)$ if $k\geq 1$
and $[g_x(0)!]:=1$. 
In the special case where the $g_x$'s are linear $\bar\nu_{V_n}^N$
is a product of Poisson probability measures conditioned on the 
hyperplane $S_N$. 

Given an energy function $H: \bbN^n \to \bbR$ we shall consider
the perturbed probability measure
\begin{equation}
\nu_{V_n}^N (\eta) = \frac{\bar\nu_{V_n}^N(\eta)}{Z_{V_n}^N} \,\nep {-H(\eta)}  \,.
\label{pert}
\end{equation}
For every $x$ we shall use the notation $\eta^{x-}$ to denote the
configuration where a particle (if there) is removed from $x$:
\begin{displaymath}
  (\eta^{x-})_y:=
  \begin{cases}
    \eta_x-1 & \text{if $\eta_x>0$ and $y=x$}\\
    \eta_y   & \text{otherwise}
  \end{cases}
\end{displaymath}
We then use $\grad_x^- f$ for the gradient $f^{x-} - f$, with $f^{x-}(\eta) =
f(\eta^{x-})$.  
The Markov generator 
\begin{equation}
\cL f = \frac1n\sum_{x,z} g_x \,\nep{-\grad_x^-H}\,\grad_{xz}f
\label{gen11}
\end{equation}
defines 
a reversible dynamics for $\nu_{V_n}^N$. Indeed, setting 
\begin{equation}
c_x(\eta) = \frac1n\,g_x(\eta_x)\,\nep{-\grad_x^-H(\eta)}\,,
\label{rateo}
\end{equation}
it is easily verified that 
the detailed balance condition holds:
\begin{equation}
  \label{eq:db1}
  \nu_{V_n}^N(\eta)c_x(\eta)=c_z(\eta^{xz})\nu_{V_n}^N(\eta^{xz})\,\quad x,z\in V_n\,.
\end{equation}
The following identity, valid for every $x,z\in V_n$ with $x\neq z$
and every function $\varphi:S_N\to\bbR$, 
is also easily verified
  \begin{equation}
    \label{eq:6}
    \nu\left[c_x\varphi\right]=
    \nu\left[c_z\varphi^{zx}\right].
  \end{equation}
Note that $c_x(\eta)=0$ iff $\eta_x=0$. 
In the dynamics defined by (\ref{gen11}) particles are removed from
$x$ with rate $g_x\,\nep{-\grad_x^-H}$ and they 
instantaneously reappear at a
uniformly chosen vertex $z\in V_n$.  
\subsection{Main estimate}
We observe that the process defined by (\ref{gen11}) can be 
written in the general frame of expression (\ref{generator}) with 
$G=\{xz\,:\;x,z\in V_n\}$, and the
rates given by  
$c(\eta,xz) =
c_x(\eta)$ for every $x,z\in V_n$. 
%
%
%
To exploit the general 
computations of the previous sections we are going to verify the following
facts.
\begin{lemma}
  \label{lemma:1}
For every $x,y,z,v \in V_n$, set
  \begin{displaymath}
    R(\eta,xz,yv):= c_x (\eta)c_y^{x-}(\eta)\nu_{V_n}^N(\eta)\,.
      \end{displaymath}
Then properties (A1)-(A4) are satisfied.
\end{lemma} 
\bpr
(A1) is trivial because both $S_N$ and $G$ are finite sets.
(A2) comes from the fact that if 
$x\neq y$ and $\eta_x\eta_y\neq0$ 
then $(\eta^{xz})^{yv}=(\eta^{yv})^{xz}$ 
while  $(\eta^{xz})^{xv}=(\eta^{xv})^{xz}$ if $\eta_x>1$.

Property (A3) holds because of the symmetry
$R(\eta,xz,yv)=R(\eta,yv,xz)$. This is obvious when $x=y$. For
$x\neq y$ it follows from 
\begin{multline*}
c_x (\eta)c_y^{x-}(\eta) = \frac1{n^2}\,
g_x(\eta_x) g_y(\eta_y)
\,\nep{-\grad_x^-H(\eta) - \grad_y^-H (\eta^{x-})} 
 = \frac1{n^2}\,
g_x(\eta_x)
g_y(\eta_y)\,
\nep{- \grad_y^- \grad_x^-H(\eta)}  \\ = c_y (\eta)c_x^{y-}(\eta) \,.
\end{multline*}
For property (A4) define, for $\eta_y>0$
\begin{displaymath}
    r_{x,y}(\eta):=
\frac{c_y(\eta^{x-})}{c_y(\eta)} \,,
  \end{displaymath}
so that $R(\eta,xz,yv)=r_{x,y}(\eta)c_x(\eta)c_y(\eta)\nu_{V_n}^N(\eta)$,
independent of $z,v$. 
Then use reversibility \eqref{eq:db1} to get 
\begin{multline*}
  \sum_{x,z,y,v}\sum_\eta R(\eta,xz,yv)F(\eta,xz,yv) =
  \sum_{x,z,y,v}\nu_{V_n}^N\left[c_x c_y r_{x,y} F(\cdot,xz,yv)\right]\\
 =\sum_{x,z,y,v}\nu_{V_n}^N\left[c_x c_y^{xz} r_{z,y}^{xz} 
F^{xz}(\cdot,zx,yv)\right].
\end{multline*}
The last term is equal to
\begin{displaymath}
  \sum_{x,z,y,v}\sum_\eta R(\eta,xz,yv) F(\eta^{xz},zx,yv)
\end{displaymath}
since it is straightforward to show that 
$c_xc_y^{xz}r_{z,y}^{xz} = c_x c_y r_{x,y}$.
\epr

Thus we can use 
Corollary~\ref{corollary1} and Corollary~\ref{corollary2} 
to bound from below the spectral gap of ${\cL}$.
We formulate the result in terms of the matrix 
\begin{equation}
\cM_{x,y}(\eta): = n\,\sqrt{c_x(\eta) c_y(\eta)} \,(1-r_{x,y}(\eta))
\label{mxy}
\end{equation}
We also use the notation
$$\varepsilon_x(\eta):=\sum_{y:\;y\neq x}|1-\nep{-\grad_x^-\grad_y^- H(\eta)}|\,.$$
\begin{theorem}\label{teom}
Assume there exists $ \d>0$ such that $\cM \geq \d$, pointwise as
quadratic forms. Then $\mathrm{gap}({\cal L})\geq \d$.
In particular, 
\begin{equation}
    \label{cobound}
\mathrm{gap}({\cal L})\geq \min_{x\in V_n}\;\min_{\eta\in
  S_N:\,\eta_x>0}
n\, c_x(\eta) \left( 1 -\frac{c_x(\eta^{x-})}{c_x(\eta)} - \varepsilon_x(\eta)
\right)\,. 
  \end{equation}
\end{theorem}
\bpr
By Lemma~\ref{lemma:1}, 
Corollary~\ref{corollary1} and 
Corollary~\ref{corollary2} we have 
\begin{displaymath}
\nu_{V_n}^N[(\cL f)^2] \geq 
\sum_{x,y,z,v}\nu_{V_n}^N\left[c_x c_y(1-r_{x,y})(\nabla_{xz}f)
(\nabla_{yv}f)\right] = \frac1n \,\nu_{V_n}^N[(u,\cM u)]\,,
\end{displaymath}
where we use the notation $(u,\cM u)=\sum_{x,y\in V_n} u_x \cM_{x,y}
u_y$, with the vectors 
$$
u_x:= \sum_{z\in V_n} \sqrt{c_x}\, \grad_{xz}f\,.
$$
By the assumption $\cM\geq \d$ we then have 
$$
 \nu_{V_n}^N[(\cL f)^2] \geq \frac\d{n}\,\nu_{V_n}^N[(u,u)]\,,
$$
where $(u,u)=\sum_x u_x^2$. To prove $\mathrm{gap}({\cal L})\geq \d$,
all we have to show is that 
\begin{equation}
\nu_{V_n}^N[(u,u)] 
= n\,\nu_{V_n}^N[f(-\cL f)]\,.
\label{medie}
\end{equation}  
This can be proved as in \cite{Ca:Po} Lemma~2.5. Namely, we rewrite
  \begin{align*}
\nu_{V_n}^N[(u,u)]  &=
    \sum_{x,z,v}\nu_{V_n}^N\left[c_x\left(\nabla_{xz}f\right)\left(\nabla_{xv}f\right)\right]\\
    & =\sum_{x,z,v}\nu_{V_n}^N\left[c_x\left(\nabla_{xz}f\right)f^{xv}\right]
    -\sum_{x,z,v}\nu_{V_n}^N\left[c_x\left(\nabla_{xz}f\right)f\right].
  \end{align*}
  The second term in the last line equals $n\,\nu_{V_n}^N[f(-\cL f)]$, while the first is 0.  In fact
  \begin{displaymath}
    \sum_{x,z,v}\nu_{V_n}^N\left[c_x\left(\nabla_{xz}f\right)f^{xv}\right]=
    \sum_{x,z,v}\nu_{V_n}^N\left[c_xf^{xz}f^{xv}\right]-
    \sum_{x,z,v}\nu_{V_n}^N\left[c_xff^{xv}\right]
  \end{displaymath}
  and by \eqref{eq:6}
  \begin{displaymath}
    \sum_{x,z,v}\nu_{V_n}^N\left[c_xf^{xz}f^{xv}\right]=
    \sum_{x,z,v}\nu_{V_n}^N\left[c_zff^{zv}\right]=
    \sum_{x,z,v}\nu_{V_n}^N\left[c_xf^{xz}f\right].
  \end{displaymath}

We turn to the proof of (\ref{cobound}). 
For any vector $w=\{w_x\}$ we have
$$(w,\cM w) = \sum_{x} n(c_x-c_x^{x-}) w_x^2 
+ \sum_x\sum_{y\neq x} n\sqrt{c_xc_y}\,(1-\nep{-\grad_x^-\grad_y^-
  H}) 
w_x w_y\,.
$$
We then estimate
$$
\sqrt{c_xc_y}\, |w_x w_y|
\leq \frac12 \left(c_x \, w_x^2 + 
c_y \,
w_y^2\right)\,.
$$
Summing over $x$ and $y\neq x$ we see that, pointwise in $\eta$: 
\begin{equation*}
(w,\cM w) \geq \sum_{x} n\, c_x \left( 1 -\frac{c_x^{x-}}{c_x} -
  \varepsilon_x\right)\, w_x^2\,,
\end{equation*}
which implies the conclusion.
\epr

\subsection{Examples}
The first observation is that when $H=0$ Theorem \ref{teom} allows
to recover exactly the result of \cite{Ca:Po} on the spectral gap of the
zero--range process under the assumption
of uniformly increasing rates. 
Indeed, if $H=0$ we have $r_{x,y}=1$ unless
$x=y$, so that $\cM$ is diagonal with entries given by 
$g_x(\eta_x) - g_x(\eta_x-1)$ which gives $\cM\geq \d$ as soon as
\begin{equation}
\min_{x\in V_n} \inf_{k\geq 0} 
\,\left[g_x(k+1)-g_x(k)\right] \geq \d\,.
\la{dd00}
\end{equation}
 
We now turn to applications of Theorem \ref{teom} 
to non--zero--range models. 
%
A class of
examples is obtained by taking the function $H$ of the form
\begin{equation}
H(\eta) = \sum_{x,y} J_{x,y} \eta_x\eta_y\,,
\label{quadr}
\end{equation}
where $J_{x,y}=J_{y,x}$ is a symmetric, constant, $n\times n$
matrix. Here $\grad_x^-H = -\sum_z J_{x,z}\eta_z$ and
$\grad_x^-\grad_y^-H = J_{x,y}$ so that
$\varepsilon_x:=\sum_{z\neq x}|1-\nep{-J_{x,z}}|\,.$
The estimate (\ref{cobound}) then becomes
\begin{equation}
\mathrm{gap}({\cal L})\geq \min_{x\in V_n}\;\min_{\eta\in S_N:\,\eta_x>0}
\,\nep{\sum_{z}J_{x,z}\eta_z}\,\left[g_x(\eta_x)-g_x(\eta_x-1)\nep{-J_{x,x}}
-g_x(\eta_x)\varepsilon_x\right]\,.
\label{jbound}
\end{equation}
\bex{zero}
The above applies in particular to the following situation.
Assume $J_{x,y}\geq 0$ for all $x,y\in V_n$. Assume also that there exists
$K\in\bbN$
such that for all $x\in V_n$ we have $J_{x,y}\neq 0$ for at most $K$
vertexes $y\neq x$. Set 
$$
a=\min_x J_{x,x}\,,\quad\; b=\max_{y\neq  x} J_{x,y}\,.
$$
Then $\varepsilon_x\leq K(1-\nep{-b})$. Assume also that we have non-decreasing
  rates:
\begin{equation}
\label{non-dec}
g_x(k+1)\geq g_x(k)\,,\quad \;k\in\bbN\,.
\end{equation}
Since 
$g_x(\eta_x)\geq 1$ and $\sum_z J_{x,z} \eta_z \geq a$ 
for any $\eta$ such that $\eta_x\geq 1$, (\ref{jbound}) gives
\begin{equation}
 \mathrm{gap}({\cal L})\geq \nep{a}\left[1 - \nep{-a}
-K(1-\nep{-b})\right]\,.
\label{jjbound}
\end{equation}
For every given $a>0$ we may take $b$ sufficiently small to obtain a
positive gap. 
\eex

\br{constantrates}
The above example includes 
the interesting case of {\em constant}
rates where 
\begin{equation}
\label{constant}
g_x(k) =  1\quad\text{ for all $k\geq 1$}\,.
\end{equation}
It is worthwhile observing that in this case if $H=0$ the gap is of
order $(1+\r)^{-2}$ with $\r=N/n$, as recently shown in
\cite{Morris} by Morris. Note that the choice (\ref{constant})
makes the reference measure $\bar\nu_{V_n}^N$ {\em uniform}
over $S_N$. Thus (\ref{jjbound}) proves that the addition of a
small mass ($a>0$) is sufficient to give a density--independent lower bound 
on the gap (for $b$ small).
\er

\bex{uno}
Here is a special case of the class of models included
in Example \ref{zero}. In particular, we assume non--decreasing rates 
as in (\ref{non-dec}). 
Consider a box of linear size $L$ in $\bbZ^d$, some $d\geq 1$, 
with periodic boundary
conditions, \ie\ we look at the quotient graph
$\Lambda_L=(\bbZ/L\bbZ)^d$. We have $n=L^d$ vertexes 
and particles jump from $x\in \Lambda$ to an arbitrary $z\in \Lambda$
with rate $c_x$ as in (\ref{rateo}) with the energy $H$ defined by
$$
H(\eta) = \b \sum_{x\sim y} \eta_x\eta_y  + \l \sum_{x}
\eta_x^2\,,
$$
where $\b,\l>0$ and the first sum runs over all pairs of adjacent
vertexes of $\Lambda$. In this case we have the expression
(\ref{quadr}) with
$$J_{x,x}= \l\,,\quad\; J_{x,y}= \begin{cases}
\frac\b{2}& \text{if} \;x\sim y\\
0 & \text{otherwise}
\end{cases}
$$ 
Since $K=2d$ here, (\ref{jjbound})
shows that 
$$
\mathrm{gap}({\cal L})\geq 
\nep{\l}\left[1-\nep{-\l} - 
2d(1-\nep{-\b})\right] \,.
$$
For every fixed $\l>0$ we can make the last expression positive by
taking $\b$ sufficiently small. 

When $\l=0$, on the other hand, (\ref{jbound}) gives useful bounds
only if we have increasing rates. Namely, set
\begin{equation}
 \e(N) := \min_x \min_{1\leq k \leq N} \frac{g_x(k)-g_x(k-1)}{g_x(k)}\,.
\label{epn}
\end{equation}
Then, if $\l=0$ (\ref{jbound}) implies 
$$\mathrm{gap}({\cal L})\geq  \left(\e(N) -
  2d(1-\nep{-\b/2})\right)\,.$$ 
This is bounded below by 
\eg\ $\e(N)/2$ as soon as $\b\leq
c\,\e(N)$ for a sufficiently small constant $c>0$. In the Poisson case
$g_x(k)=k$, $\e(N)=1/N$ so that $\b$ has to be taken as small as
  $O(1/N)$. 
Clearly, if the rates grow exponentially, 
\eg\ $g_x(k)=\nep{k}$ we have $\e(N)$ bounded away from zero
  independently
of $N$ (this is like having a mass again). 
 \eex


\section{Kawasaki-type dynamics in the continuum}
\label{kaw:cont}
In this section we consider a system of particles jumping
about a bounded subset of $\bbR^d$. In many respects the model
described below may be considered as the continuous version of
the random walk models of Section 3 and Section 4.

Let $\O$ be the set of locally finite subsets of $\R^d$. We provide $\O$ with the weakest topology that, for every continuous $f :\R^d \ra \R$ with compact support, makes the maps $\eta \mapsto \sum_{x \in \eta} f(x)$ continuous. Measurability on $\O$ is provided by the corresponding Borel $\s$-field.

Now let $\L$ be a bounded Borel subset of $\R^d$ of nonzero Lebesgue measure, and set
\[
S := \O_{\L} := \{\eta \in \O: \eta \subseteq \L\}.
\]
Consider a nonnegative measurable and even function $\varphi:\R^d \ra [0,+\infty)$ (everything works with minor modifications for $\varphi:\R^d \ra [0,+\infty]$ allowing ``hardcore repulsion''). We fix a {\em boundary condition} $\tau \in \O_{\L^c} := \{ \eta \in \O: \eta \subseteq \L^c\}$, and define the Hamiltonian $H^{\tau}_{\L}:S \ra [0,+\infty]$
\begin{equation}
H^{\tau}_{\L}(\eta) = \sum_{\stackrel{\{x,y\} \subseteq \eta \cup \tau}{\scriptscriptstyle \{x,y\} \cap \L \neq \emptyset}} \varphi(x-y).
\label{hamcont}
\end{equation}
The dependence of $H^{\tau}_{\L}$ on $\L$ and $\tau$ is omitted in the sequel.

For $N\in\natural$ we let $S_N=\{\eta \in S : |\eta| = N\}$
denote the subset of $S$ consisting of all possible configurations
of $N$ particles in $\Lambda$. Note that a measurable function $f:S_N \ra \R$ may be identified with a symmetric function from $\L^N \ra \R$. With this identification, we assume that the boundary condition $\tau$ is such that $H(\eta) < +\infty$ in a subset of $\L^N$ having positive Lebesgue measure.
Now, for $\b>0$, we define the {\em canonical} Gibbs measure in the finite volume $\L$ with inverse temperature $\b$ as the probability $\nu_{\L}^N$ on $S_N$ given by
\begin{displaymath}
  \nu_\Lambda^N[f]:=
  \frac{1}{Z_\Lambda^N}\int_{\Lambda^N}\frac{dw}{|\Lambda|^N}\;e^{-\b H(w)}f(w)\,,
\end{displaymath}
for any bounded function $f:S_N\to \bbR$, where $Z_{\Lambda}^N$ is a normalization factor.

For $x,z\in\Lambda$ define the map on $S$:
\begin{displaymath}
  \gamma_{xz}(\eta):=
  \begin{cases}
    \eta\setminus\{x\}\cup\{z\} &\text{if $x\in\eta$} \\
    \eta & \text{otherwise.}
  \end{cases}
\end{displaymath}
Define the map
$\gamma_{x}^-(\eta)=\eta\setminus\{x\}$ (if $x\in\eta$, otherwise
$\gamma_{x}^-(\eta)=\eta$). 

As usual we set $G:=\{\gamma_{xz}:x,z\in\Lambda\}$. 
In the sequel we will write 
$\eta^{xz}$ for $\gamma_{xz}(\eta)$,  
$\eta^{x-}$ for $\gamma_{x-}(\eta)$, $\nabla_{xz}$ for 
$\nabla_{\gamma_{xz}}$, and $\nabla_{x}^-$ for $\nabla_{\gamma_{x}^-}$.
Furthermore, for any function $f$ on $S$ we define 
$f^{xz}(\eta):=f(\eta^{xz})$ and $f^{x-}(\eta):=f(\eta^{x-})$.

Consider the following Markov generator
\begin{displaymath}
  ({\cal L}f)(\eta):=
  \sum_{x\in\eta}\int_\Lambda \frac{dz}{|\Lambda|}\,\,e^{-\b(H^{xz}(\eta)-H^{x-}(\eta))}\,\nabla_{xz}f(\eta).
\end{displaymath}
In words, this corresponds to moving particles $x\in \eta$ \,to 
a point $z\in \L$ with infinitesimal rate $$\frac1{|\L|}\,\,e^{-\b(H^{xz}(\eta)-H^{x-}(\eta))}\,dz\,.$$ 
Observe that 
\begin{equation}
H^{xz}(\eta)-H^{x-}(\eta) = \sum_{y\in\eta\setminus \{x\}} \varphi(y-z)\,.
\label{varphio}
\end{equation}
It can be shown that ${\cal L}$ has a domain of self-adjointness in 
$L^2(\nu_\Lambda^N)$, and that generates a Markov semigroup.
The core ${\cal{C}}$ can be taken as the set of bounded functions 
$f:S_N\to \bbR$.

This generator is of the form (\ref{generator}) if we define 
$c(\eta,d\g)$ by
\[
\int c(\eta,d\g)F(\g) := \sum_{x \in
  \eta}\int_{\L}\frac{dz}{|\Lambda|}\,\,
e^{-\b(H^{xz}-H^{x-})} F(\g_{xz}).
\]
In particular, it is easy to show that the reversibility 
condition (\ref{db}) holds.
The Dirichlet form associated with ${\cal L}$ is
\begin{equation}
  {\cal E}(f,f)=
  \frac{1}{2}\int_\Lambda \frac{dz}{|\Lambda|}\,\,\nu_\Lambda^N\left[\sum_{x\in\eta}e^{-\b(H^{xz}-H^{x-})}(\nabla_{xz}f)^2\right].
\label{dirfor}
\end{equation}
\begin{lemma}
  \label{lemma:2}
  Define
  \begin{displaymath}
    r(\cdot,xz,yv):=
    \begin{cases}
      0 & \text{if $x=y$} \\
      e^{-\b\,\varphi(z-v)} & \text{if $x\neq y$},
    \end{cases}
  \end{displaymath}
and $R(\cdot,\gamma,\delta):=\nu_\Lambda^N r(\cdot,\gamma,\delta)c(\cdot,d\gamma)c(\cdot,d\delta)$.
Then (A1)-(A4) are satisfied.
\end{lemma}
 
\bpr
Property (A1) is a consequence of the
fact that $R$ is bounded (recall that $\varphi\geq 0$). Therefore any
bounded function is in $L^1(R)$. 
(A2) comes from the fact that if $x\neq y$ then $\nu_\Lambda^N$-almost surely $(\eta^{xz})^{yv}=(\eta^{yv})^{xz}$.
(A3) holds because $r(\eta,xz,yv)=r(\eta,yv,xz)$.
Property (A4) can be checked as follows:
\begin{multline}
\label{eq:3}
  \int\nu_\Lambda^N(d\eta) r(\eta,\gamma,\delta)c(\eta,d\gamma)c(\eta,d\delta)F(\gamma(\eta),\delta,\gamma)=\\
  =\int_\Lambda\frac{dz}{|\Lambda|}\int_\Lambda\frac{dv}{|\Lambda|}\,
\nu_\Lambda^N\left[\sum_{x\in\eta}\sum_{y\in\eta\setminus\{x\}}
e^{-\b [H^{xz}-H^{x-}+H^{yv}-H^{y-}+\varphi(z-v)]}
F(\eta^{xz},\gamma_{xz},\gamma_{yv})\right]=\\
=\frac{N(N-1)}{Z_\Lambda^N|\Lambda|^{N+2}}
\int_{\Lambda^N}dw
\int_{\Lambda}dz
\int_{\Lambda}dv\,\, e^{-\b[H^{w_1
  z}-H^{w_1-}+H^{w_2v}-H^{w_2-}+H+
\varphi(z-v)]}F^{w_1z}(\cdot,\gamma_{w_1z},\gamma_{w_2v})\,,
\end{multline}
where \eg
$$
H^{w_1 z}(w)=H^{w_1 z}(w_1,\dots,w_N):=H(z,w_2,\dots,w_N)$$
and $H^{w_2-}:=H(w_1,w_3,w_4,\dots,w_N)$.
By the change of variables $w_1\mapsto z\mapsto w_1$ we see that the last term in \eqref{eq:3} equals
\begin{displaymath}
 \frac{N(N-1)}{Z_\Lambda^N|\Lambda|^{N+2}}\int_{\Lambda^N}
dw\int_{\Lambda}dz\int_{\Lambda}dv\; 
e^{-\b[H-H^{w_1-}+H^{w_1z,w_2v}-H^{w_1z,w_2-}+H^{w_1z}+\varphi(w_1-v)]}F(\cdot,zw_1,w_2v)\,.
\end{displaymath}
Since
$$
  H^{w_1z,w_2v}-H^{w_1z,w_2-}+\varphi(w_1-v)=
  H^{w_2v}-H^{w_2-}+\varphi(v-z)\,,
$$
the above implies (A4).
\epr
We define two parameters
\begin{equation}
\e_1 = 
\sup_{\eta\in S_{N-1}} \int_{\L}\frac{dv}{|\L|}\,\,
\left(1 - e^{-\b\sum_{x\in\eta}\varphi(v-x)}\right)\,,\quad\;
\e_2 = 
2\,(N-1)\,\sup_{z\in\L}\int_{\L}\frac{dv}{|\L|}\,\,
\left(1 - e^{-\b\,\varphi(v-z)}\right)\,.
\label{epss}
\end{equation}
\bt{teo_con_cont}
 For any non--negative $\varphi$, $\L\subset\bbR^d$ a bounded
 Borel set and $N\in\bbN$, $\beta\geq 0$ we have 
\begin{displaymath}
    \gap({\cal L})\geq
    1\,-\,\e_1\,-\,\e_2\,.
  \end{displaymath}
  \et
\bpr
By Lemma~\ref{lemma:1}, Corollary~\ref{corollary1} and Corollary~\ref{corollary2} we have to bound from below
\begin{displaymath}
  \int_{\Lambda}\frac{dz}{|\Lambda|}
\int_{\Lambda}\frac{dv}{|\Lambda|}
\,\nu_\Lambda^N\left[\sum_{x,y\in\eta}(1-r(\cdot,xz,yv))\,
e^{-\b(H^{xz}-H^{x-}+H^{yv}-H^{y-})}(\nabla_{xz}f)(\nabla_{yv}f)\right]\,,
\end{displaymath}
in terms of the Dirichlet form $\cE(f,f)$.
The above can be written as $A+B$ where
\begin{multline}
\label{eq:400}
A = \int_{\Lambda}\frac{dz}{|\Lambda|}
\int_{\Lambda}\frac{dv}{|\Lambda|}
\,\nu_\Lambda^N\left[\sum_{x\in\eta}\,
e^{-\b(H^{xz}-H^{x-} +H^{xv}-H^{x-})}(\nabla_{xz}f)(\nabla_{xv}f)\right]\,\\
= \frac{N}{Z_\Lambda^N|\Lambda|^{N+2}}\int_{\Lambda^N}
dw\int_{\Lambda}dz\int_{\Lambda}dv\; 
e^{-\b(H^{w_1z}-H^{w_1-}+H^{w_1v}-H^{w_1-}+H)}(\nabla_{w_1z} f)
(\nabla_{w_1v} f)\,,
\end{multline}
and
\begin{multline}
\label{eq:401}
B= \int_{\Lambda}\frac{dz}{|\Lambda|}
\int_{\Lambda}\frac{dv}{|\Lambda|}
\,\nu_\Lambda^N\left[\sum_{x,y\in\eta:\;y\neq x}(1-r(\cdot,xz,yv))\,
e^{-\b(H^{xz}-H^{x-}+H^{yv}-H^{y-})}(\nabla_{xz}f)(\nabla_{yv}f)\right]
\\ = \frac{N(N-1)}{Z_\Lambda^N|\Lambda|^{N+2}}
\int_{\Lambda^N}dw\int_{\Lambda}dz\int_{\Lambda}dv\; 
e^{-\b(H^{w_1 z}-H^{w_1-}+H^{w_2v}-H^{w_2-}+H)}
\left[1-e^{-\b\,\varphi(z-v)}\right](\nabla_{w_1z} f)(\nabla_{w_2v} f)\,.
\end{multline}
We next show that 
\begin{equation}
  \label{eq:5}
  A=
    \frac{N}{2Z_\Lambda^N|\Lambda|^{N+1}}\int_{\Lambda^N}dw
\int_{\Lambda}dz\; e^{-\b(H^{w_1 z}-H^{w_1-}+H)}(\nabla_{w_1z} f)^2\int_{\Lambda}\frac{dv}{|\Lambda|}e^{-\b(H^{w_1v}-H^{w_1-})}.
\end{equation}
In fact, 
using a change of variables as in Lemma
\ref{lemma:2}
we see that 
 \begin{align*}
- \int_{\Lambda^N} &dw\int_{\Lambda}dz\int_{\Lambda}dv\; 
e^{-\b(H^{w_1 z}-H^{w_1-}+H^{w_1v}-H^{w_1-}+H)}(\nabla_{w_1z} f)f
\\ &\quad\quad
=\frac12\int_{\Lambda^N}dw\int_{\Lambda}dz\; 
e^{-\b(H^{w_1 z}-H^{w_1-}+H)}(\nabla_{w_1z} f)^2
\int_{\Lambda}dv\,e^{-\b(H^{w_1v}-H^{w_1-})},
 \end{align*}
and
\begin{multline*}
\int_{\Lambda^N}dw\int_{\Lambda}dz\int_{\Lambda}dv\; 
e^{-\b[H^{w_1 z}-H^{w_1-}+H^{w_1v}-H^{w_1-}+H]}(\nabla_{w_1z} f)f^{w_1v}=\\
  =\int_{\Lambda^N}dw\int_{\Lambda}dz\int_{\Lambda}dv\; 
e^{-\b(H^{w_1 z}-H^{w_1-}+H^{w_1v}-H^{w_1-}+H)}f^{w_1z}f^{w_1v}+\\
  -\int_{\Lambda^N}dw\int_{\Lambda}dz\int_{\Lambda}dv\; 
e^{-\b(H^{w_1 z}-H^{w_1-}+H^{w_1v}-H^{w_1-}+H)}ff^{w_1v}=0.
\end{multline*}
Decomposing $\nabla_{w_1v} f = f^{w_1v} - f$ in (\ref{eq:400})
this proves (\ref{eq:5}). Recalling (\ref{dirfor}) and
(\ref{varphio}) we then see that 
\begin{equation}
\label{402}  
A\geq \inf_{w\in\Lambda^N}
\left[\int_{\Lambda}\frac{dv}{|\Lambda|}\,
e^{-\b(H^{w_1v}-H^{w_1-})}\right]\,{\cal E}(f,f)
\geq (1- \e_1)\,{\cal E}(f,f)\,.
\end{equation}

We now estimate the absolute value of $B$ in (\ref{eq:401}) from above. Using
$$|\nabla_{w_1z} f \nabla_{w_2v} f|\leq \frac12\,[(\nabla_{w_1z} f)^2
+(\nabla_{w_2v} f)^2]$$ and $e^{-\b(H^{w_2v}-H^{w_2-})}\leq 1$ we
easily obtain
\begin{align*}
  B & \geq -\frac{N(N-1)}{Z_\Lambda^N|\Lambda|^{N+1}}\int_{\Lambda^N}
dw\int_{\Lambda}dz\; e^{-\b(H^{w_1 z}-H^{w_1-}+H)}
(\nabla_{w_1z} f)^2
\int_{\Lambda}\frac{dv}{|\Lambda|}\left[1-e^{-\beta\,\varphi(z-v)}\right]
\\
& \geq -\e_2\,{\cal E}(f,f)
\end{align*}
Together with (\ref{402}) this completes the proof of the theorem.
\epr

Similarly to what will be seen in the non--conservative case treated in section \ref{sec:gl}
an application of the above results shows that a positive gap
is obtained under high--temperature/small--density assumptions.
 
We first observe that for fixed $\L$ and $N$ we have $\e_1,\e_2\to 0$
as $\b\to 0$, so that $\gap({\cal  L})\to 1$ by Theorem \ref{teo_con_cont}. 
To obtain quantitative estimates involving the density
of particles $N/|\L|$ we may use the following criterion.

\bc{corre}
Assume that the non--negative pair potential
$\varphi$ and the inverse temperature $\b$ satisfy
$$
\e(\b):=\int_{\bbR^d}\left(1-e^{-\b\varphi(x)}\right)\,dx \,<\,\infty\,.
$$
Then, for every bounded Borel set $\L\subset\bbR^d$, $N\in\bbN$
$$
\gap({\cal  L}) \geq 1\,-\,\frac{3(N-1)}{|\L|}\,\e(\b)\,.
$$
\ec
\bpr
Let $\e_1,\e_2$ be as in Theorem \ref{teo_con_cont}.
Clearly, $$\e_2\leq \frac{2(N-1)}{|\L|}\,\int_{\bbR^d}\left(1-e^{-\b\varphi(x)}\right)\,dx\,.$$ Moreover,
using the elementary inequality
$$
1-e^{-s-t}\leq (1-e^{-s}) + (1-e^{-t})\,,\quad s,t\geq 0\,,
$$
we see that $$\e_1\leq \frac{N-1}{|\L|}\,\int_{\bbR^d}\left(1-e^{-\b\varphi(x)}\right)\,dx\,.$$
\epr

\section{Glauber dynamics with unbounded, discrete spin} \label{GD}

In this section we consider a multidimensional birth and death process. Given a finite set $\L$ (no geometrical structure is required for the moment), we let $S := \N^{\L}$. Thus, for $\eta = (\eta_x)_{x \in \L}$, $\eta_x$ denotes the number of particles at the site $x \in \L$. We consider the creation an annihilation maps on $S$: for $x \in \L$
\[
[\g^+_x(\eta)]_y = \begin{cases} \eta_x +1 & \mbox{for } y=x \\ \eta_y & \mbox{otherwise}
\end{cases}
\]
\[
[\g^-_x(\eta)]_y = \begin{cases} \eta_x -1 & \mbox{if } y=x \mbox{ and } \eta_x >0 \\ \eta_y & \mbox{otherwise.}
\end{cases}
\]
We let $G:= \{\g^+_x,\g^-_x: x \in \L\}$. In the sequel we write $\nabla^+_x$ and $\nabla^-_x$ rather than $\nabla_{\g^+_x}$ and $\nabla_{\g_x^-}$. We consider a birth and death process with generator of the form
\be{gengla}
{\cal L}f(\eta) := \sum_{x \in \L} \left[c(\eta,\g_x^+) \nabla^+_x f(\eta) + c(\eta,\g_x^-) \nabla^-_x f(\eta)\right],
\ee
where $c(\eta,\g_x^+)$ is the rate of creation of a particle at $x$, and $c(\eta,\g_x^-)$ is the rate of annihilation of a particle at $x$. Let $\nu$ be a probability on $S$ such that $\nu(\eta)>0$ for every $\eta \in S$. We set
\begin{eqnarray}
c(\eta,\g_x^+) & := & (\eta_x + 1) \frac{\nu(\g_x^+ \eta)}{\nu(\eta)} \label{rate+gla} \\ c(\eta,\g_x^-) & := & \eta_x \label{rate-gla}.
\end{eqnarray}
With these rates we have that $\left(\g_x^+\right)^{-1} = \g_x^-$, $\left(\g_x^-\right)^{-1} = \g_x^+$ in the sense of condition (Rev) (although the equality $\g_x^+ (\g_x^- \eta)$ fails if $\eta_x=0$).
Moreover the detailed balance condition
\[
c(\eta,\g_x^+) \nu(\eta) = c(\g_x^+ \eta,\g_x^-) \nu(\g_x^+ \eta)
\]
holds, which is equivalent to (\ref{db}) for this case. The measure $R$ is chosen according to Proposition \ref{proposition1}, with $J := \{\g_x^+; x \in \L\}$. Note that $J \cap J^{-1} = \emptyset$. More explicitly:
\begin{eqnarray}
r(\eta,\g_x^+,\g_y^+) &  = & \frac{c(\g^+_x \eta, \g^+_y)}{c(\eta, \g^+_y)} \nonumber\\
r(\eta,\g_x^-,\g_y^-) &  = & \frac{c(\g^-_x \eta, \g^-_y)}{c(\eta, \g^-_y)} = \begin{cases} 0 & \mbox{if } \eta_x \eta_y =0 \\ 1 & \mbox{if } \eta_x \eta_y \neq 0, x \neq y \\ \frac{\eta_x -1}{\eta_x} & \mbox{if } \eta_x \neq 0, x=y \end{cases}
 \label{rgla} \\
r(\eta,\g_x^-,\g_y^+) & = & r(\eta,\g_x^+,\g_y^-) \ = \ 1. \nonumber
\end{eqnarray}
Note that 
\[
{\cal D}({\cal L}) = \left\{ f \in L^2(\nu) : \nu\left[c(\eta,\g_x^+)\left( \nabla^+_x f(\eta)\right)^2 \right] < +\infty \ \forall x \in \L\right\}.
\]
As a core ${\cal{C}}$ for ${\cal L}$ we take
\be{coregla}
{\cal{C}}:= \left\{ f \in L^2(\nu): \exists N>0 \mbox{ such that } \nabla^+_xf (\eta) = 0 \ \forall x \in \L \mbox{ for } \sum_{x \in \L} \eta_x >N\right\}.
\ee
\bl{glal1}
For the measure $R$ given in (\ref{rgla}), properties (A1)-(A4) are satisfied.
\el
\bpr
For the above choice of ${\cal{C}}$, for $f \in {\cal{C}}$ the map $(\eta,\g,\d) \mapsto \nabla_{\g}f(\eta) \nabla_{\d}f(\eta)$ has a bounded support, so (A1) is easily satisfied. Properties (A2) and (A4) follow from Proposition \ref{proposition1}. Property (A3) comes from the fact that $r(\eta,\g,\d) = r(\eta,\d,\g)$ for every $(\eta,\g,\d) \in S \times G \times G$, as is easily checked from (\ref{rate+gla}), (\ref{rate-gla}) and (\ref{rgla}).
\epr

\subsection{Example: pair interaction in a Poissonian field}

We assume here $\nu$ is of the following form:
\[
\nu(\eta) := \frac{1}{Z} \prod_{x \in \L} \frac{\l^{\eta_x}}{\eta_x !} \exp\left[- \b\sum_{\{x,y\} \cap \L \neq \emptyset} \varphi(x,y,\eta_x, \eta_y) \right],
\]
where $\L$ is a finite subset of $\Z^d$, $\b>0$ and
\[
\varphi: \Z^d \times \Z^d \times \N \times \N \ra \R
\]
is a pair potential, such that $\varphi(x,y,m, n) = \varphi(y,x,n,m)$ for every $x,y \in \Z^d$, $n,m \in \N$, and $\varphi(x,x,n,m)  \equiv 0$. The measure $\nu$ on $S$ depend on the {\em boundary condition} $\eta\big|_{\L^c}$, that is supposed to be equal to a given fixed $\tau \in \N^{\L^c}$; this dependence is omitted in the notation. For the above measure to be well defined for every choice of boundary condition we require that, for every $x \in \L$, $\eta \in \N^{\Z^d}$, the infinite sum
\[
\sum_{y \in \Z^d} \varphi(x,y,\eta_x, \eta_y)
\]
is well defined and takes value in $(-\infty, +\infty]$. For example, this holds true in either one of the following cases:
\bi
\item
$\varphi$ is nonnegative;
\item
$\varphi$ is of finite range, \ie there exists $k>0$ such that $\varphi(x,y,\eta_x, \eta_y) \equiv 0$ for $|x-y| >k$.
\ei
With this choice of $\nu$ the rates become
\[
c(\eta,\g_x^-) = \eta_x, \hspace{2cm} c(\eta,\g_x^+) = \l \exp\left[ -\b\sum_{y \in \Z^d }  \nabla_x^+
\varphi(x,y,\eta_x,\eta_y) \right].
\]
\bt{gla}
Define
\[
\e(\b) := \sup_{\stackrel{x \in \Z^d}{\scriptscriptstyle \eta \in \N^{\Z^d}}}  \left\{ \sum_{z \in \Z^d} \exp\left[ -\b\sum_{y \in \Z^d }  \nabla_z^+
\varphi(z,y,\eta_z,\eta_y) \right]\left| 1- \exp\left[ -\b \nabla_x^+
\nabla_z^+  \varphi(z,x,\eta_z,\eta_x) \right]  \right|  \right\}.
\]
Then
\[
\gap({\cal L}) \geq1-\l \e(\b).
\]
Note that this bound is independent of $\L$ and of the boundary condition.
\et
\bpr
We first observe that, since $r(\eta,\g,\d) =1$ for $(\g,\d) \in \left(J \times J^{-1}\right) \cup \left(J^{-1} \times J \right)$,
\begin{multline*}
\int \nu(d\eta)c(\eta,d\g)
c(\eta,d\d)[1-r(\eta,\g,\d)] \nabla_{\g}f(\eta)\nabla_{\d}f(\eta) \\ = \int_{J \times J} \nu(d\eta)c(\eta,d\g)
c(\eta,d\d)[1-r(\eta,\g,\d)] \nabla_{\g}f(\eta)\nabla_{\d}f(\eta) \\ + \int_{J^{-1} \times J^{-1}} \nu(d\eta)c(\eta,d\g)
c(\eta,d\d)[1-r(\eta,\g,\d)] \nabla_{\g}f(\eta)\nabla_{\d}f(\eta).
\end{multline*}
By (\ref{rgla})
\begin{multline*}
\int_{J^{-1} \times J^{-1}} \nu(d\eta)c(\eta,d\g)
c(\eta,d\d)[1-r(\eta,\g,\d)] \nabla_{\g}f(\eta)\nabla_{\d}f(\eta) \\ = \sum_{x \in \L} \nu\left[\eta_x \left(\nabla_x^- f(\eta)\right)^2 \right] =  \E(f,f),
\end{multline*}
where we have used (\ref{dirj}).
On the other hand
\begin{multline*}
\int_{J \times J} \nu(d\eta)c(\eta,d\g)
c(\eta,d\d)[1-r(\eta,\g,\d)] \nabla_{\g}f(\eta)\nabla_{\d}f(\eta) \\
= \sum_{x \in \L} \nu\left[c(\eta, \g^+_x) \sum_{z \in \L} c(\eta,\g^+_z) \left( 1- \exp\left[ -\b \nabla_x^+
\nabla_z^+ \sum_y \varphi(z,y,\eta_z,\eta_y) \right] \right) \nabla_x^+f(\eta)  \nabla_z^+f(\eta)\right] \\
= \sum_{x \in \L} \nu\left[c(\eta, \g^+_x) \sum_{z \in \L} c(\eta,\g^+_z) \left( 1- \exp\left[ -\b \nabla_x^+
\nabla_z^+  \varphi(z,x,\eta_z,\eta_x) \right] \right) \nabla_x^+f(\eta)  \nabla_z^+ f(\eta)\right] 
\end{multline*}
Thus, by Schwarz inequality, 
\[
\left|\int_{J \times J} \nu(d\eta)c(\eta,d\g)
c(\eta,d\d)[1-r(\eta,\g,\d)] \nabla_{\g}f(\eta)\nabla_{\d}f(\eta)\right| \leq C,
\]
where
\[
C = \sum_{x \in \L} \nu\left[c(\eta, \g^+_x) \sum_{z \in \L} c(\eta,\g^+_z) \left| 1- \exp\left[ -\b \nabla_x^+
\nabla_z^+  \varphi(z,x,\eta_z,\eta_x) \right]  \right| \left(\nabla^+_xf(\eta)  \right)^2\right].
\]
Now, the inequality
\[
C \leq \l \e(\b) \E(f,f),
\]
is rather immediate from the definition of $\e(\b)$ and (\ref{dirj}), and thus
\[
\int \nu(d\eta)c(\eta,d\g)
c(\eta,d\d)[1-r(\eta,\g,\d)] \nabla_{\g}f(\eta)\nabla_{\d}f(\eta) \geq  [1-\l\e(\b)] \E(f,f).
\]
The conclusion now follows from Corollary \ref{corollary2}.
\epr

The lower bound on the spectral gap given in Theorem \ref{gla}, depends on the inverse temperature $\b$ and on the density $\l$ of the reference Poissonian field. We now give an example where the estimate on $\e(\b)$ can be carried out explicitly.
\bex{exgla1}
Let $K: \Z^d \ra [0,+\infty)$ be such that $K(0) =0$, $K(-x) = K(x)$ and
\[
\sum_{x \in \Z^d} K(x) < +\infty,
\]
and define
\[
\varphi(x,y,\eta_x,\eta_y) := K(x-y)\eta_x \eta_y .
\]
This example is consistent with the interpretation of a configuration $\eta \in \N^{\Z^d}$ as a system of particles in $\Z^d$: each pair of particles gives a positive contribution $K(x-y)$ to the interaction energy, that depends on the relative position $x-y$ of particles. Since adding one particle increases the interaction energy, $\nabla_z^+
\varphi(z,y,\eta_z,\eta_y) \geq 0$, and therefore
\[
 \exp\left[ -\b\sum_{y \in \Z^d }  \nabla_z^+
\varphi(z,y,\eta_z,\eta_y) \right] \leq 1.
\]
Moreover $\nabla_x^+
\nabla_z^+  \varphi(z,x,\eta_z,\eta_x) = K(x-z)$, so that
\[
\e(\b) \leq \sum_{z \in \Z^d} \left(1- e^{-\b K(z)}\right).
\]
In particular we have that $\e(\b) \ra 0$ as $\b \ra 0$. Thus, the condition $\l \e(\b) <1$, which guarantees a positive spectral gap, is a high temperature and/or small density condition, as one would expect.

\eex

\section{Glauber dynamics of particles in the continuum}
\label{sec:gl}
As we mentioned in the introduction, for the models we describe in this section, estimates for the spectral gap were obtained via the Bakry-Emery approach in \cite{Ko:Ly}. Our aim here is to show that this computation falls within our general scheme. 

We use here the same notations introduced in Section \ref{kaw:cont}. In addition, we assume the nonnegative pair potential $\varphi$ and the inverse temperature $\b$ to satisfy 
 the condition
\be{integrability}
\e(\b):= \int_{\R^d} \left(1-e^{-\b \varphi(x)}\right)dx < +\infty.
\ee

Functions from $S$ to $\R$ may be identified with symmetric functions from $\bigcup_n \L^n$ to $\R$. With this identification, we define the finite volume {\em grand canonical} Gibbs measure $\nu_{\L}$ with inverse temperature $\b>0$ and activity $z>0$ by
\begin{equation}
\nu_{\L}[f] := \frac{1}{Z} \sum_{n=0}^{+\infty} \frac{z^n}{n!} \int_{\L^n}
e^{-\b H(x)}f(x)dx,
\label{nuf}
\end{equation}
where $Z$ is the normalization.

As in Section \ref{GD} we define the creation an annihilation maps on $S$: for $x \in \L$
\[
\g^+_x(\eta) = \eta \cup \{x\}
\]
\[
\g^-_x(\eta) = \eta \setminus \{x\}.
\]
We let $G:= \{\g^+_x,\g^-_x: x \in \L\}$. In the sequel we write $\nabla^+_x$ and $\nabla^-_x$ rather than $\nabla_{\g^+_x}$ and $\nabla_{\g_x^-}$. Note that $\nabla^-_x f(\eta) = 0$ unless $x \in \eta$. We consider the following Markov generator
\be{genglac}
{\cal L}f(\eta) := \sum_{x \in \eta} \nabla^-_xf(\eta) + z \int_{\L} e^{-\b\nabla^+_x H(\eta)} \nabla^+_x f(\eta).
\ee
It is shown in \cite{Be:Ca:Ce}, Proposition 2.1, that ${\cal L}$ has a domain of self-adjointness in $L^2(\nu_{\L})$, and that generates a Markov semigroup. It is also shown that a core ${\cal{C}}$ is given by
\be{coreglac}
{\cal{C}}:= \{ f \in L^2(\nu_{\L}): \exists \, M>0 \mbox{ such that } |f| \leq M \mbox{ and } f(\eta) = 0 \mbox{ for } |\eta| > M\},
\ee
where $|\eta|$ denote the cardinality of $\eta$. This generator is indeed of the form (\ref{generator}) if we define $c(\eta,d\g)$ by
\[
\int F(\g)c(\eta,d\g) := \sum_{x \in \eta} F(\g^-_x) + z \int_{\L} e^{-\b\nabla^+_x H(\eta)} F(\g^+_x)dx.
\]
In particular, it is easy to show that the reversibility condition (\ref{db}) holds.

Similarly to Section \ref{GD}, the measure $R$ is chosen according to Proposition \ref{proposition1}, with $J := \{\g_x^+; x \in \L\}$. In particular
\begin{eqnarray}
r(\eta,\g^+_x,\g^+_y) & = & \frac{dc(\g^+_x \eta, \cdot)}{dc( \eta, \cdot)}(\g^+_y) \ = \ \exp\left[- \b \nabla^+_x \nabla^+_y H(\eta)\right] \ = \ \exp\left[-\b\varphi(x-y)\right] \nonumber \\
r(\eta,\g^-_x,\g^-_y) & = & \frac{dc(\g^-_x \eta, \cdot)}{dc( \eta, \cdot)}(\g^-_y)  \ = \ \begin{cases} 1 & \mbox{for } x,y \in \eta, x \neq y \\ 0 & \mbox{otherwise} \end{cases} \label{rglac} \\
r(\eta,\g_x^-,\g_y^+) & = & r(\eta,\g_x^+,\g_y^-) \ = \ 1. \nonumber
\end{eqnarray}
\bl{glacl1}
For the measure $R$ defined in (\ref{rglac}) properties (A1)-(A4) hold.
\el
\bpr
For property (A1), note that the function $r(\eta,\g,\d)$ in (\ref{rglac}) is bounded. Therefore it is enough to prove that, for $f \in {\cal{C}}$, the function $(\eta,\g,\d) \mapsto \nabla_{\g}f(\eta) \nabla_{\d}f(\eta)$ is in $L^1(\nu_{\L}(d\eta)c(\eta,d\g)c(\eta,d\d)$. But
\begin{multline*}
\int \left|\nabla_{\g}f(\eta) \nabla_{\d}f(\eta)\right| \nu_{\L}(d\eta)c(\eta,d\g)c(\eta,d\d) = \nu_{\L}\left[\left(\int \left|\nabla_{\g}f(\eta)\right|c(\eta,d\g)\right)^2\right] \\
=\nu_{\L}\left[\left(\sum_{x \in \eta} \left|\nabla^-_x f(\eta)\right| + z\int_{\L} e^{-\b\nabla^+_x H(\eta)} \left| \nabla^+_x f(\eta) \right| \right)^2 \right].
\end{multline*}
The last integrand
\[
\left(\sum_{x \in \eta} \left|\nabla^-_x f(\eta)\right| + z\int_{\L} e^{-\b\nabla^+_x H(\eta)} \left| \nabla^+_x f(\eta) \right| \right)^2
\]
is bounded, since $f$ is bounded and supported on sets up to a certain cardinality, and, by non negativity of the potential $\varphi$, $e^{-\b\nabla^+_x H(\eta)} \leq 1$. This completes the proof for property (A1). Properties (A2) and (A4) follow from Proposition \ref{proposition1}, while (A3) comes from the symmetry property $r(\eta,\g,\d) = r(\eta,\d,\g)$.
\epr
\bt{glac}
Let $\e(\b)$ be the quantity defined in (\ref{integrability}). Then
\[
\gap({\cal L}) \geq [1-z\e(\b)].
\]
Note that this bound is independent of $\L$ and the boundary condition $\tau$.
\et
\bpr
The proof is quite close to the one of Theorem \ref{gla}. We begin observing that
\begin{multline*}
\int \nu_{\L}(d\eta)c(\eta,d\g)
c(\eta,d\d)[1-r(\eta,\g,\d)] \nabla_{\g}f(\eta)\nabla_{\d}f(\eta) \\ = \int_{J \times J} \nu_{\L}(d\eta)c(\eta,d\g)
c(\eta,d\d)[1-r(\eta,\g,\d)] \nabla_{\g}f(\eta)\nabla_{\d}f(\eta) \\ + \int_{J^{-1} \times J^{-1}} \nu_{\L}(d\eta)c(\eta,d\g)
c(\eta,d\d)[1-r(\eta,\g,\d)] \nabla_{\g}f(\eta)\nabla_{\d}f(\eta).
\end{multline*}
By (\ref{rgla})
\begin{multline*}
\int_{J^{-1} \times J^{-1}} \nu_{\L}(d\eta)c(\eta,d\g)
c(\eta,d\d)[1-r(\eta,\g,\d)] \nabla_{\g}f(\eta)\nabla_{\d}f(\eta) \\
 =
\nu_{\L}\left[ \sum_{x \in \eta} \left(\nabla^-_x f(\eta)\right)^2 \right] = \E(f,f),
\end{multline*}
where we have used (\ref{dirj}). On the other hand
\begin{multline*}
\int_{J \times J} \nu_{\L}(d\eta)c(\eta,d\g)
c(\eta,d\d)[1-r(\eta,\g,\d)] \nabla_{\g}f(\eta)\nabla_{\d}f(\eta) \\ = \nu_{\L}\left[z^2 \int_{\L^2} e^{-\b \nabla^+_x H(\eta)} e^{-\b \nabla^+_y H(\eta)} \left(1-e^{-\b \varphi(x-y)}\right) \nabla^+_x f(\eta) \nabla^+_y f(\eta) dxdy \right].
\end{multline*}
Thus, by Schwarz inequality
\[
\left|\int_{J \times J} \nu_{\L}(d\eta)c(\eta,d\g)
c(\eta,d\d)[1-r(\eta,\g,\d)] \nabla_{\g}f(\eta)\nabla_{\d}f(\eta) \right| \leq C,
\]
where
\begin{multline*}
C \leq  \nu_{\L}\left[z^2 \int_{\L} e^{-\b \nabla^+_x H(\eta)}\left[ \int_{\L}e^{-\b \nabla^+_y H(\eta)} \left(1-e^{-\b \varphi(x-y)}\right)dy\right]\left( \nabla^+_x f(\eta)\right)^2  dx \right] \\ \leq z\e(\b)  \nu_{\L}\left[z \int_{\L} e^{-\b \nabla^+_x H(\eta)}\left( \nabla^+_x f(\eta)\right)^2  dx \right] ,
\end{multline*}
where we used (\ref{dirj}) and the fact that $e^{-\b \nabla^+_y H(\eta)} \leq 1$. The conclusion now follows readily as in Theorem \ref{gla}.
\epr

\noindent
{\em Acknowledgments}. We thank Prof. M. Ledoux for several discussions concerning this work. We also thank the Associate Editor for useful comments and suggestions


\begin{thebibliography}{99}

\bibitem{Ba:Em}
Bakry, D. and Emery, M., {\em Diffusions hypercontractives}. S\'eminaire de Probabilit\'es XIX. Lecture Notes in Math. 1123, 177-206, Springer-Verlag, 1985.

\bibitem{Be:Ca:Ce}
Bertini, L., Cancrini, N. and Cesi, F., {\em The spectral gap for a Glauber-type dynamics in a continuous gas}, Ann. I. H. Poincar\'e - PR 38, 1 (2002), 91-108.

\bibitem{Bo}
Bochner, S., {\em Vector fields and Ricci Curvature}, Bull. Amer. Math. Soc., 52 (1946), 776-797.

\bibitem{Ca:Ma}
Cancrini, N. and Martinelli, F., {\em On the spectral gap of Kawasaki dynamics under a mixing condition revisited}, J. Math. Phys. 41 (2000), n. 3, 1391-1423.

\bibitem{Ca:Ma:Ro}
Cancrini, N., Martinelli, F.  and Roberto, C.,
{\em The logarithmic Sobolev constant of Kawasaki dynamics under a mixing condition revisited},
Ann. Inst. H. Poincar{\'e} Probab. Statist.  38  (2002),  no. 4, 385--436.

\bibitem{Ca:Po}
Caputo, P., Posta, G., {\em Entropy dissipation estimates in a zero
  range dynamics}, preprint {\tt math.PR/0405455}.

\bibitem{Ca:St}
Carlen E. A.,  Stroock D. W., {\em An application of the Bakry-\'Emery criterion to infinite dimensional diffusions}, Sem. de Probabilit\'es XX, Springer-Verlag, LNMS 1204 (1986), 341-348.

\bibitem{Da:Pa:Po}
Dai Pra, P., Paganoni A.M. and Posta G.,
{\em Entropy inequalities for unbounded spin systems},
The Annals of Probability, Vol. 30 no 4 (2002), 1959-1976.

\bibitem{De:St}
Deuschel, J.-D. and Stroock, D. W.,
{\em Hypercontractivity and spectral gap of symmetric diffusions with
applications to the stochastic Ising models}.
 J. Funct. Anal.  92  (1990),  no. 1, 30--48.


\bibitem{He}
Helffer, B., {\em Remarks on Decay of Correlations and Witten Laplacians
Brascamp Lieb Inequalities and Semiclassical Limit}, J. Funct. Anal. 155 (1998), 571-586.

\bibitem{Ko:Ly}
Kondratiev, Y. and Lytvynov, E., {\em Glauber Dynamics of continuous
  particle systems} preprint {\tt math.PR/0306252}

\bibitem{Le}
Ledoux, M., {\em Logarithmic Sobolev inequalities for unbounded spin systems revisited}, S\'eminaire de Probabilit\'es XXXV. Lecture Notes in Math. 1755, 167-194, Springer, Berlin,  2001.

\bibitem{Lic}
Lichn\'erowicz, A., {\em G\'eom\'etrie des groupes de transformations}, Dunod, Paris, 1958.

\bibitem{Li}
Liggett T.M.,  {\em Interacting particle systems}.
276. Springer-Verlag, New York-Berlin, 1985.

\bibitem{Lu:Ya}
Lu, S. T. and Yau, H.-T., 
{\em Spectral gap and logarithmic Sobolev inequality for Kawasaki and Glauber dynamics},
Communications in Mathematical Physics 156 (1993), 399--433

\bibitem{ma:ol}
Martinelli, F. and Olivieri, E., {\em Approach to equilibrium of Glauber dynamics in
the one phase region. II. The general case},  Comm. Math. Phys. 161 (1994), 487-514.

\bibitem{Morris} Morris, B., 
{\em Spectral gap for the zero range process with constant rate}, preprint
{\tt math.PR/0405161}

\bibitem{St}
Stroock, D. W., {\em An Introduction to the Analysis of Paths on a Riemannian Manifold}, American Mathematical Society, 1999.

\bibitem{sz1}
Stroock, D. W. and Zegarlinski B. , {\em The logarithmic Sobolev inequality for
discrete spin systems on a lattice}, Comm. Math. Phys. 149 (1992), 175-193.
\bibitem{sz2}
Stroock, D. W. and Zegarlinski B. , {\em  The logarithmic Sobolev inequality for
continuous spin systems on a lattice}, J. Funct. Anal. 104 (1992),  299-326.

\bibitem{Ya}
Yau, H-T., {\em Logarithmic Sobolev Inequality for Lattice Gases with Mixing Conditions}, Comm. Math. Phys. 181 (1996), 367-408.


\end{thebibliography}
\end{document}